\documentclass{article}
\usepackage{arxiv}

\usepackage{graphicx}
\usepackage{amssymb}

\usepackage{lineno}

\usepackage{amsmath}
\usepackage{hyperref}

\providecommand{\myfloor}[1]{\left \lfloor #1 \right \rfloor }
\usepackage{amsfonts}
\usepackage{epstopdf}
\usepackage{float} 
\usepackage{caption} 
\usepackage{subcaption} 
\usepackage{enumitem} 

\usepackage{amssymb}
\usepackage{amsthm} 
\usepackage{amsmath}

\newenvironment{rcases}
  {\left.\begin{aligned} \displaystyle}
  {\end{aligned}\qquad \right\rbrace}

\title{A Fast, Spectrally Accurate Homotopy based Numerical Method for Solving Nonlinear Differential Equations}

\author{
  Andrew C. Cullen
    \\
  School of Mathematical Sciences\\
  Monash University\\
  Australia \\
  \texttt{andrew.cullen@monash.edu} \\
   \And
Simon R. Clarke \\
  School of Mathematical Sciences\\
  Monash University\\
  Australia\\
   \\
}

\begin{document}
\maketitle

\begin{abstract}
We present an algorithm for constructing numerical solutions to one--dimensional nonlinear, variable coefficient boundary value problems. This scheme is based upon applying the Homotopy Analysis Method (HAM) to decompose a nonlinear differential equation into a series of linear differential equations that can be solved using a sparse, spectrally accurate Gegenbauer discretisation. Uniquely for nonlinear methods, our scheme involves constructing a single, sparse matrix operator that is repeatedly solved in order to solve the full nonlinear problem. As such, the resulting scheme scales quasi--linearly with respect to the grid resolution. We demonstrate the accuracy, and computational scaling of this method by examining a fourth--order nonlinear variable coefficient boundary value problem by comparing the scheme to Newton--Iteration and the Spectral Homotopy Analysis Method, which is the most commonly used implementation of the HAM.\\
\end{abstract}


\section{Introduction}

The Homotopy Analysis Method (HAM) \citep{Liao1992} is a modern technique for iteratively constructing analytic solutions to nonlinear equations, that can be considered as an extension of the \cite{Lyapunov1892} and Homotopy Perturbation methods \citep{He2003}. Applying this scheme to vexatious nonlinear problems has demonstrated that the HAM has advantageous convergence properties, relative to those exhibited by other techniques \citep{Liao1997, Liao1998, Liao2002b}. The technique works by decomposing a nonlinear problem into a sequentially coupled series of linear differential equations, through the use of a homotopy, which is a concept rooted in differential geometry. However, as an analytic method, it relies on being able to solve the linear equations either by hand, or through the use of computer algebra packages. This requirement limits the scope of problems that can be approached using the HAM. As such, there is particular interest in developing new numerical schemes, that leverage the convergence properties of the HAM. As a consequence of the method's convergence properties, a particular focus has been placed upon high-order accurate methods, which can produce fundamentally more accurate and reliable solutions for numerical domains discretised at lower resolutions, which in turn impacts upon the computational cost of these schemes.

These high-order accurate methods can in turn be broadly categorised into local and global (or spectral) methods, with the later involving discretisations of derivatives that depend upon all the points in the numerical domain, while local methods calculate derivatives at each point in the domain in terms of its adjacent elements \citep{Strikwerda2004}. Local high-order methods are employed in multi-block finite difference methods, some finite volume methods, and stabilised finite element methods; whereas Discontinuous Glarkerin (DG), DG spectral element methods, and spectral volume and difference methods all fall under the categorisation of global methods \citep{Chan2017, Wang2002, Yibiao2016}. These schemas are typically considered as the basis for solving unsteady problems upon a wide range of domains. However, while these schemes can be applied to solve steady nonlinear problems by coupling them within an iterative framework, the schemes can introduce significant computational costs. Furthermore, these schemes are not guaranteed to remain stable for sufficiently long to resolve the steady solution.

While a number of specific schemes to approach steady, nonlinear problems have been proposed, including Newton iteration,  AiTEM \citep{Yang2007}, nonlinear Conjugate Gradient \citep{Dhai1999} and Multigrid based methods \citep{Hemker1990, Ronquist1987}, Newton iteration remains amongst the most commonly employed. These schemes all can be considered as linear solvers that are iterated upon in order to converge upon the motivating nonlinear problem. However, as was the case with the earlier examples, the iterative schemes that result often scale poorly with the grid resolution.

This paper we will outline a new approach to constructing numerical solutions through the HAM, called the Gegenbauer Homotopy Analysis Method (GHAM), and examine its numerical performance relative to other HAM based approaches, and other numerical schemes for solving nonlinear differential equations. We first describe the HAM in Section 2, in the context of its original formulation as a semi--analytic scheme, and as a numerical scheme when solved using the Spectral Homotopy Analysis Method (SHAM). In Section 3,  we introduce the GHAM, and discuss its advantages over previously developed nonlinear solvers. Section 4 will present a theoretical justification of the difference in numerical performance between the GHAM and other techniques. Finally, in Section 5 we show numerical results, convergence, and scaling properties for a fourth--order nonlinear variable coefficient boundary value problem, and compare these results for the Gegenbauer- and Spectral-Homotopy Analysis methods with Newton Iteration and MATLAB's inbuilt boundary value problem routine `BVP4C'. These comparisons will allow us to understand the potential for the GHAM to significantly reduce the computational cost of solving steady, nonlinear boundary value problems. 

\section{Homotopy Analysis Method}\label{ch:homotopy}

The Homotopy Analysis Method (HAM) is a technique for solving nonlinear equations by constructing a homotopy that defines a smooth, continuous deformation from one equation onto another. If a homotopy can be found between the two equations, it follows that the same homotopy should also hold for the solutions of the two equations, and thus allowing one known solution to an equation to be deformed onto the unknown solution of another equation. Within the context of the HAM, when searching for a solution to a motivating nonlinear problem, a second equation can be introduced, for which there is a known solution, and then this solution can be deformed onto the motivating nonlinear problem through the homotopy.  

If we consider two continuous functions $f$ and $g$ which are respectively the products of the topological spaces $X$ and $Y$ then the homotopy is defined as being a continuous function $\mathcal{H}: X \times [0,1] \to Y$, such that if $x \in X$ then $\mathcal{H}(\textbf{x},0) = f(\textbf{x})$ and $\mathcal{H}(\textbf{x},1) = g(\textbf{x})$. For the purposes of solving a nonlinear differential equation, $f(\textbf{x})$ and $g(\textbf{x})$ can be thought to be the approximate and exact solutions respectively. 

The process of constructing this homotopy involves partitioning a nonlinear problem into an infinite sequence of linear sub--problems, which can then be solved sequentially. One crucial feature is that the technique introduces very few limitations upon the linear sub--problems being employed in the solution process. Introducing an additional convergence control parameter allows the technique to be applied to problems that have either been considered intractable, or that have been limited to narrow regions of convergence. From a fluid mechanics perspective these enhanced convergence control properties have been used to further explore Blasius boundary layers \citep{Liao2002}, Magnetohydrodynamic flow \citep{Hayat2007}, steady-state resonant progressive waves \citep{Xu2012} and Von Karman swirling flow in the presence of viscosity \citep{Yang2006}.

Considering this technique in terms of the nonlinear differential operator $\mathcal{N}$ where

\begin{equation}\label{eqn:originalN}
\mathcal{N} \left[ u(\textbf{x}) \right] = \psi(x),
\end{equation}

the solution to this can be constructed by deforming from a solution of an arbitrary linear differential equation $\mathcal{L}[u(\textbf{x})]$ through the homotopy
\begin{equation}\label{eqn:originalhomotopy}
\mathcal{H} \lbrack \phi(x);q \rbrack = (1-q) \mathcal{L} \lbrack \phi(\textbf{x};q) - u_{0}(\textbf{x};q) \rbrack - q \hbar \left( \mathcal{N} \lbrack \phi(\textbf{x};q) \rbrack - \psi(x) \right), 
\end{equation}
A specific feature of the HAM is the introduction of the auxiliary convergence control parameters $\hbar$. The homotopy itself is controlled by the deformation parameter $q \in [0,1]$ and $\phi(\textbf{x};q)$ is a representation of the solution across $q$, corresponding to the deformation from a trial function $u_{0}(\textbf{x};q)$ to $u(\text{x})$. The function $\phi(\textbf{x};q)$ can be considered as a Maclaurin series in terms of $q$ 
\begin{equation}
\label{eqn:psi}
\begin{rcases}
\phi(\textbf{x};q) &= \phi(\textbf{x};0) + \sum_{m=1}^{\infty} \frac{1}{m!} \frac{\partial^{m} \phi(\textbf{x};q)}{\partial q^m} \bigg|_{q=0} q^{m}\\
& = U_{0} + \sum_{m=1}^{\infty} U_{m}(\textbf{x}) q^{m}, 
\end{rcases}
\end{equation}
where 
$$U_{m}(\textbf{x}) = \frac{1}{m!} \frac{\partial^{m} \phi(\textbf{x};q)}{\partial q^m} \bigg|_{q=0}.  $$
By imposing without loss of generality that $\mathcal{H} = 0$, the behaviour of $\phi$ can be considered further in the context of the limits of $q$ within Equation (\ref{eqn:originalhomotopy}). Therefore 
\begin{equation}\label{eqn:basichomotopy}
\begin{rcases}
&\mathcal{H}[\phi(\textbf{x};0);0] = \mathcal{L}[\phi(\textbf{x};0) - u_{0}(\textbf{x})],\qquad && \phi(\textbf{x};0) = u_{0}(\textbf{x}),\\
&\mathcal{H} \lbrack \phi(\textbf{x};1);1 \rbrack = \mathcal{N} \lbrack \phi(\textbf{x};1) \rbrack - \psi(x),\qquad && \phi(\textbf{x};1) = u(\textbf{x}).
\end{rcases}
\end{equation}
This corresponds to deforming $\phi(\textbf{x};q)$ from an initial guess of $u_{0}(\textbf{x})$ onto the exact solution $u(\textbf{x})$. This is contingent upon $\phi(\textbf{x};q)$ being analytic at $q=0$, and more broadly existing over $q \in [0,1]$.
The existence and convergence of the scheme outlined in Equation (\ref{eqn:originalhomotopy}) is determined by the choice of the convergence control parameter $\hbar$, and the auxiliary linear operator $\mathcal{L}$.  The latter parameter can be defined almost arbitrarily, although our testing has indicated that the most efficient choice of $\mathcal{L}$ resembles the linearisation of $\mathcal{N}[u(\textbf{x})]$---an idea that will be explored later.

Subject to the choice of $\mathcal{L}$, the homotopy formulation of Equation (\ref{eqn:originalhomotopy}) can be employed to solve an equation of the form of (\ref{eqn:originalN}) by substituting $\phi$ from Equation (\ref{eqn:psi}) into the homotopy equation (\ref{eqn:originalhomotopy}), differentiating $m$-times with respect to $q$ and then setting $q = 0$. This then results in the sequence of serially dependent linear differential equations
\begin{equation}
\label{eq:Homotop}
\begin{rcases}
& \mathcal{L} \lbrack U_{m} (\textbf{x}) - \chi_{m} U_{m-1}(\textbf{x}) \rbrack = \hbar H(\textbf{x}) R_{m},\\
&R_{m} = \frac{1}{(m-1)!} \left\lbrace  \frac{\partial^{m-1}}{\partial q^{m-1}} \mathcal{N} \left[ \sum_{n=0}^{m-1} U_{n} (\textbf{x})q^{n} \right] \right\rbrace \bigg|_{q=0} - (1-\chi_{m})\psi(x), \\
&\chi_{m} = \begin{cases}
0 & \text{if $m \leq 1$,}\\
1 & \text{if $m \geq 2$.}
\end{cases}
\end{rcases}
\end{equation}
The sum within $R_{m}$ has been truncated to $n \in [0,m-1]$ as the terms of order $q^{m}$ and higher will all vanish after setting $q = 0$, irrespective of the form of the nonlinearity. As a consequence of this, Equation (\ref{eq:Homotop}) must be strictly linear with respect to $U_{m}$. This result means that through careful choice of $\mathcal{L}$, the system can be iteratively solved through analytic, semi--analytic or numerical schemes.

As the homotopy parameter $q$ maps the original equation at $q = 0$ to the motivating nonlinear equation at $q=1$, it follows that $u(\textbf{x}) = \phi(\textbf{x}; 1)$ in Equation (\ref{eqn:psi}), i.e.,
\begin{equation}\label{eqn:u_form}
u(\textbf{x}) = U_{0} + \sum_{m=1}^{\infty} U_{m}(\textbf{x}).
\end{equation}

\subsection{Spectral Homotopy Analysis Method}\label{sec:SHAM}

While the original formulation of the HAM was built upon solving Equation (\ref{eq:Homotop}) algebraically, the scheme has since evolved to incorporate both computer algebra and numerical discretisations. The latter approach using spectral methods falls under the aegis of the Spectral Homotopy Analysis Method (SHAM). These semi--analytic and numerical extensions were driven by one of the fundamental limitations of the HAM---in that the auxiliary linear operator $\mathcal{L}$ must be chosen in such a way that the resulting equations can be solved. Computer algebra approaches broaden the forms of $\mathcal{L}$ that can be solved for, with the numerical framework increasing the set of admitted forms of $\mathcal{L}$ even further, to the point where it is no longer a primary consideration when solving a nonlinear equation.

\cite{Motsa2010, Motsa2012} was one of the first authors to consider strictly numerical analogues of the HAM, using Chebyshev collocation matrices to develop what is now known as the SHAM. Initial testing of this scheme demonstrated that the SHAM has the potential to significantly outperform MATLAB's inbuilt boundary value problem routine `BVP4C' \citep{Nik2013}, although a formal exploration of its numerical scaling has not been conducted.

When discretised numerically, the iterative process of the SHAM can be expressed as
\begin{equation} \label{eqn:SHAMNumerical}
\textbf{U}_{m} = \chi_{m}\textbf{U}_{m-1} + \hbar \textbf{A}^{-1} \left[\textbf{N}_{m}\right].
\end{equation}
Here $\textbf{A}$ is the discretised matrix operators for the linear differential equation using the Chebyshev collocation approach, and $\textbf{N}_{m}$ is the matrix representation of $R_{m}$. Subject to a choice of $\hbar$, the general numerical solution to the original nonlinear problem is then
\begin{equation}
\textbf{U} = \sum_{m=0}^{\infty} \textbf{U}_{m}.
\end{equation}
This approach to numerically discretise the linear differential equations that result from using the HAM shares many of the same advantages and disadvantages of other collocation based schema. While it does exhibit spectral accuracy in solving the resulting linear problems, it does so at the cost of having to solve dense matrix systems. Beyond this, the scheme is unable to handle variable coefficient differential equations, since the coefficients exhibit large condition numbers and become almost singular. This final limitation is a particular problem when domain mappings from semi--infinite or infinite domains are incorporated into the problem, in order to solve over the Chebyshev domain $x \in [-1,1]$.\\ 

\section{Gegenbauer Homotopy Analysis Method}\label{subsection:GHAM}

For all their advantageous numerical properties, Chebyshev polynomials have limited utility for constructing matrix operators for variable coefficient linear boundary value problems, as the resulting matrix operators rapidly become singular with increases in the grid resolution. An alternative approach is to turn to the Gegenbauer (otherwise known as Ultraspherical) polynomials, of which the Chebyshev polynomials are a subset. These Gegenbauer polynomials exhibit the previously introduced advantageous properties of the Chebyshev polynomials, while also yielding sparse matrix operators.

Similar to the Chebyshev polynomials, the Gegenbauer polynomials have an orthogonality property, in this case with respect to the weight $(1 - x^2)^{\lambda - \frac{1}{2}},$ and take the form
$$C_{k}^{(\lambda)}(x) = \frac{2^{k} (\lambda + k - 1)!}{k! (\lambda - 1)!} x^{k} + \mathcal{O}(x^{k-1}).$$
More specifically, the Gegenbauer polynomials can be constructed in terms of the recurrence relationship 
\begin{equation}\label{eqn:gegenbauerrecurrence}
C_{j+1}^{(\lambda)} (x) = \frac{2 (j + \lambda)}{j + 1} x C_{j}^{(\lambda)}(x) - \frac{j + 2 \lambda - 1}{j+1} C_{j-1}^{(\lambda)} (x), \qquad j \geq 1,
\end{equation}
subject to $C_{0}^{(\lambda)} = 1$ and $C_{1}^{(\lambda)} = 2 \lambda x$ for all $\lambda \in \mathbb{Z}^{+}$.  One notable advantage of constructing a numerical scheme in terms of the Gegenenbauer polynomials is that differentiation in Chebyshev space can be expressed as 
\[
 \frac{d^{k} T_{j}}{d x^{k}} =
  \begin{cases} 
      2^{k-1} j(k-1)! C_{n-k}^{(k)}  \hfill & \text{ for $j \geq k$,} \\
      0 \hfill & \text{ for $0 \leq j \leq k-1$, } \\
  \end{cases}
\]
where $T_j(x)$ is the $j$-th Chebyshev polynomial. Consequently Chebyshev differentiation operators of any order can be represented by a sparse matrix in Gegenbauer space. This compares favourably to the traditional approach of constructing collocation based Chebyshev differentiation matrix operators, which results in dense matrix representations.

The sparse numerical discretisation of a generalised variable coefficient boundary value problem
\begin{equation}\label{eqn:GegenbauerLinearGeneral}
\begin{rcases}
&\mathcal{L} u(x) = f(x) \; \text{ and } \mathbf{\mathcal{B}}u(x) = \mathbf{c} \\
&\mathcal{L} = a_{N} (x) \frac{d^{N}}{dx^{N}} + \cdots + a_{1}(x) \frac{d}{dx} + a_{0}(x)
\end{rcases}
\end{equation}
can be achieved by following \cite{Olver2013}, with the numerical operator taking the form 
\begin{equation}\label{eqn:assembledGeg}
\mathcal{L} := \mathcal{M}_{N} [ a^{N}] \mathcal{D}_{N} + \sum_{\lambda = 1}^{N-1} \mathcal{S}_{N-1} \ldots \mathcal{S}_{\lambda} \mathcal{M}_{\lambda}[a^{\lambda}] \mathcal{D}_{\lambda} + \mathcal{S}_{N-1} \ldots \mathcal{S}_{0} \mathcal{M}_{0} [a^{0}], 
\end{equation}
so that
\begin{equation}
A 
\begin{pmatrix}  
\hat{u}_{0} \\ \hat{u}_{1} \\ \vdots \\ \hat{u}_{n-1}  \end{pmatrix} 
= 
\begin{pmatrix}  
\\ \mathbf{c} \\ \\ \mathcal{P}_{n-K}\mathcal{S}_{n-1}\mathcal{S}_{n-2}\cdots \mathcal{S}_{0} \mathbf{f} \\ \hspace{1 cm}
\end{pmatrix}
\end{equation}
where 
\begin{equation}
A = \begin{pmatrix}
\mathbf{\mathcal{B}} \mathcal{P}_{n}^{T} \\
\\
\mathcal{P}_{n-K} \mathcal{L}
\end{pmatrix}.
\end{equation}
Here $\mathcal{L}$ is the matrix defined by Equation (\ref{eqn:assembledGeg}), $\mathbf{\mathcal{B}}$ represents the boundary conditions through the method of boundary bordering, $\mathcal{P}_j \in \mathbb{R}^{n \times n}$ is a projection operator $(I_{j}, \mathbf{0})$ constructed in terms of the identity matrix $I_{j}$, which has size $\mathbb{R}^{j \times j}$, and $K$ is the number of boundary conditions. The resultant matrix operators from this system are sparse, almost--banded, can be solved in $\mathcal{O}(m^2n)$ operations, and require $\mathcal{O}(mn)$ storage, where $n$ is the number of Chebyshev modes required to resolve the solution, and $m = 1 + (\text{number of non-zero superdiagonals}) + (\text{number of non-zero subdiagonals})$. The number of non--zero diagonals excludes the $K$ dense rows corresponding to the boundary conditions. 

This linear discretisation can in turn be leveraged to serve as the numerical solver for the linear sub--problems that result from applying the HAM. The nonlinear discretisation for the GHAM replicates the approach taken for the SHAM, with the difference being in the application of the boundary conditions. To understand this, let us consider a nonlinear problem of the form
\begin{equation}\label{eqn:GHAMGeneralNonlinear}
 \left.\begin{aligned}
       &\mathcal{N} \lbrack u(x) \rbrack = \psi(x), x \in (a,b)\\
		& \mathcal{B} \lbrack u(a),u'(a),\ldots,u(b),u'(b),\ldots \rbrack = \textbf{c}.
       \end{aligned}
 \right\}
\end{equation}
Inhomogeneous boundary conditions have the potential to significantly complicate the iterative process that will result from applying the HAM to this nonlinear problem. As such, homogeneous boundary conditions can be imposed by solving the problem 
\begin{equation}
\begin{rcases}
&\mathcal{L} \lbrack u_{0}(x) \rbrack = \psi(x), \\
&\mathcal{B} \lbrack u(x),u'(x),\ldots \rbrack = \textbf{c}, \\
\end{rcases}
\end{equation}
for an arbitrary auxiliary linear operator $\mathcal{L}$, to construct a mapping for Equations (\ref{eqn:GHAMGeneralNonlinear}) whereby  
$$u(x) = v(x) + u_{0}(x).$$
An alternate approach is to arbitrarily construct a $u_{0}(x)$ that satisfies the inhomogeneous boundary conditions, and apply the mapping outlined above. In either case, this substitution gives rise to the modified nonlinear equation 
$$\mathcal{L}_{1}[v(x)] + \mathcal{N}_{1}[v(x)] = \psi_{1}(x).$$
The process of decomposing $u(x)$ into an equation with homogeneous boundary conditions is not strictly necessary, however it does greatly simplify the process of enforcing the boundary conditions for each step of the homotopy process. 

Within the framework of this modified nonlinear equation, a homotopy can be constructed between the auxiliary linear operator $\mathcal{L}$ and the new nonlinear equation in the same manner as described earlier, resulting in the iterative scheme 
\begin{equation}
\textbf{V}_{m} = \chi_{m} \textbf{V}_{m-1} + \hbar \textbf{A}^{-1} \left( \textbf{A}_{1} \textbf{V}_{m-1}  + \textbf{N}_{m} \right).
\end{equation}
Within this framework, $\textbf{V}_{m}$ is the vector representation of $v_{m}$, $\textbf{A}$ and $\textbf{A}_{1}$ are the constant matrix operators corresponding to $\mathcal{L}$ and $\mathcal{L}_{1}$ respectively as represented through the Gegenbauer method, and $\textbf{N}_{m}$ represents $R_{m}$. From this, the solution to the motivating nonlinear equation (\ref{eqn:GHAMGeneralNonlinear}) is then
\begin{equation}
\textbf{U} = \textbf{U}_{0} + \sum_{m=0}^{\infty} \textbf{V}_{m}.
\end{equation}
The difference between this framework and that expressed by Equation (\ref{eqn:SHAMNumerical}), beyond the manner in which the matrices are discretised, stems from the mapping to homogeneous boundary conditions.

\section{Theoretical scaling}\label{sec:scaling}

The difference between the GHAM and other comparable numerical schemes can be considered by examining the cost of solving the matrix equations that these schemes produce. Typically, nonlinear solvers require the construction of a new dense matrix inverse at each step of the iterative process, however the Spectral- and Gegenbauer Homotopy Analysis Methods can be considered in terms of a single matrix operator, that is either dense in the case of SHAM, or sparse for GHAM. Through the Coppersmith--Winograd algorithm, the single matrix inversion required for methods based upon the HAM can be constructed with $\mathcal{O}(n^{2.373})$ operations for a matrix system in $\mathbb{R}^{n \times n}$ \citep{Davie2013,LeGall2014}. Each step of the homotopy iterative process can then be solved with an additional $\mathcal{O}(n^2)$ operations. However, taking this approach results in a dense inverse matrix, and as such the matrix--vector products are in terms of dense systems---negating the fundamental advantages of the Gegenbauer discretisation. These is also a secondary cost with respect to the amount of memory required to store these dense matrices as $n$ increases.

An alternate approach for the GHAM involves  using LU decomposition for the matrix inversion \citep{Horn1985}. A test on a fourth--order boundary value problem discretised using the Gegenbauer method over $512$ points resulted in a sparse, diagonally dominant matrix operator for which only $2.5\%$ of the matrix components are non--zero. In contrast, the LU decomposition involves a significant increase in the fill rate of both matrices, with non--zero elements covering $49 \%$ of the equivalent dense matrix.

This increase in the overall sparsity decreases the memory requirements for storing these matrices, however employing LU decomposition is still desirable for its implications in terms of the time required to solve these systems. For a dense matrix system, the cost of constructing the original LU decomposition is $\mathcal{O}(n^3)$ operations, with solutions of a matrix equation involving the LU decomposition able to be solved in $\mathcal{O}(n^2)$ operations. For dense matrices these scaling properties make LU decomposition a poor alternative to the Coppersmith--Winograd algorithm, however in the case of sparse matrices the partial preservation of sparsity within LU decomposition significantly improves the scaling of both the decomposition process, and the solution of the resulting matrix equations. For systems where the number of non--zero points is $\mathcal{O}(n)$ \cite{Gines1998} showed that the cost of both constructing the LU decomposition and evaluating the matrix-vector products reduces to $\mathcal{O}(n)$, with results to this effect being further outlined within the proceeding section.

For the systems that result from using the GHAM, it was found that the most computationally efficient and effective form of LU decomposition involves scaled pivoting---whereby the relative magnitudes of points, rather than their absolute magnitudes, are taken into account in the selection of the pivots. For the GHAM, the specific implementation of this process was handled through MATLAB's inbuilt sparse LU factoriser, which is in turn built upon the Suitesparse UMFPACK package \citep{Davis2006, Davis2011}.

Of course, solving the matrix systems is only one component of the overall computational cost of solving a nonlinear equation (or system of equations) within the GHAM framework. Other factors that influence the exhibited scaling include setting up the Gegenbauer matrix system; evaluating derivatives; the transform between Chebyshev and real space; and evaluating the right hand side term $R_{m}$. The scaling of $R_m$ is particularly important to understand, due to its scaling with $m$.

While the computational complexity for evaluating $R_{m}$ is heavily problem dependent, an understanding of the scheme's scaling can be obtained by examining the cases of strictly quadratic and cubic nonlinearities. For a problem exhibiting a strictly quadratic nonlinearity solved over $M$ steps of the GHAM iterative process, evaluating each $R_{m}$ for $m \leq M$ involves $\myfloor{\frac{m^2}{4}}$ multiplication operations and $\myfloor{\frac{(m-2)^2}{4}}$ addition operations, where $\myfloor{\cdot}$ represents the floor function. For a problem discretised over $n$ quadrature points, the computational complexity of evaluating all the $R_{m}$ terms will scale with $\mathcal{O}\left( n m^2 \right)$.

In the case of a strictly cubic nonlinearity, the number of multiplication operations required is $\mathcal{O} \left(m + \myfloor{\frac{3 (m-1)^2 + 1}{4}}\right)$, and the number of addition operations is $(m-2)^2$, and as such will scale in a similar manner to that for a strictly quadratic nonlinearity. 

For nonlinearities exhibiting more complex structures, if $G^{(r)}(U_{0}) \neq 0$ for all $r$, then evaluating each $R_{m}$ using Faa di Bruno's formula will involve an additional $h(m)$ products and $p(m)$ sums, where $p(m)$ is the number of integer partitions of $m$ and $h(m)$ is the total number of parts in all partitions of $m$. While neither of these two terms have an analytic description, they do respectively correspond to the $A000041$ and $A006128$ sequences from \cite{OEIS}. Our testing shows that the quadratic component of this appears to dominate, and that the scheme still broadly scales as $\mathcal{O}(n m^2)$.

As such, the dominant factor for solving a system of equations using the GHAM is the cost of solving the matrix equations, which, ignoring sparsity will behave as 
\begin{equation}\label{eqn:generalscalingproperties}\mathcal{O}\left(n m^2 + n^{3} + m n^{2} \right),
\end{equation}
where again $m$ is the number of iterations and $n$ is the number of grid points, and it can be assumed that $m \ll n$.

However, if we instead impose that the number of non--zero elements of the matrix equations is $\mathcal{O}(n)$ then as was described previously the LU decomposition also scales with $\mathcal{O}(n)$, which leads to Equation (\ref{eqn:generalscalingproperties}) reducing to
\begin{equation}\label{eqn:scaling_result}
\mathcal{O}(nm^2).
\end{equation}
Of course these results are for calculations at a single value of $\hbar$, and any exploration of the convergence properties across a range of the convergence control parameter will necessitate an increase in the computational cost. However the number of samples in $\hbar$ should be small relative to $n$ and $m$, and as such (\ref{eqn:scaling_result}) will still hold.

Subject to the assumption that $m$ is independent of $n$, this process demonstrates quasi-linear $\mathcal{O}(n)$ scaling with respect to the grid resolution. Within the field of numerical analysis of nonlinear differential equations, all currently employed techniques scale with, absolute minimum, $\mathcal{O}(n^2)$. These theoretical results will now be confirmed by considering the solution of a fourth--order boundary value problem.

\section{Numerical Testing of the GHAM}\label{sec:TwoDViscous}

To test the validity of the above theoretical scaling properties, we turn to the problem of a two--dimensional flow of a laminar, viscous, incompressible fluid confined within a rectangular domain bounded by moving porous walls, which can be recast as a nonlinear, variable coefficient boundary value problem. This problem has applications to mixing processes, as well as for boundary layer control systems, and was chosen as it has previously been studied by  \cite{Motsa2014} through the use of SHAM. The form of the problem considered is
\begin{equation}\label{eq:TwoDViscous}
\begin{rcases}
&y^{(iv)} + \alpha \left(\frac{x + 1}{4} y''' + \frac{3}{4} y''\right) + R_{e} \left(\frac{1}{2} y y''' - \frac{1}{4} y' y'' \right) = 0 \text{ for } x \in [-1,1] \\
&y(-1) = 0, \qquad y''(-1) = 0,\\
&y(1) = 1, \qquad \hspace{0.34 cm} y'(1) = 0.
\end{rcases}
\end{equation}
Here $\alpha$ is the non-dimensional wall dilation rate; $R_{e}$ is the permeation Reynolds number, which is positive for fluid injection across the boundary, and negative in the presence of suction on the boundaries; and $y(x)$ is a dimensionless function describing the mass flow as a function of the distance in the wall normal direction, from the reference point at $x = -1$.

As was discussed previously, the solution to a nonlinear differential equation using Homotopy based techniques can be constructed in terms of a range of auxiliary linear operators, which will serve as the basis of the iterative scheme. To explore the implications of this freedom, a range of auxiliary linear operators can be constructed and compared, which replicate the original nonlinear differential equation to varying degrees. For this problem the following operators will be used
\begin{equation}
\begin{rcases}
\mathcal{L}_{1} [y] &= y^{(iv)}, \\
\mathcal{L}_{2} [y] &= y^{(iv)} + \alpha \left(\frac{x + 1}{4} y''' + \frac{3}{4} y''\right), \\
\mathcal{L}_{3} [y] &= y^{(iv)} + \alpha \left(\frac{x + 1}{4} y''' + \frac{3}{4} y''\right) + R_{e} \left(\frac{1}{2} y''' - \frac{1}{4} y'' \right), \\
\mathcal{L}_{4} [y] &= y^{(iv)} + \alpha \left(\frac{x + 1}{4} y''' + \frac{3}{4} y''\right) + \frac{1}{2} R_{e} \left( \hat{y}_{0} y''' - \hat{y}_{0}' y'' \right).
\end{rcases}
\label{eq:aux_linear_operators}
\end{equation}
\begin{figure}
\centering
\begin{subfigure}{.5\textwidth}
  \centering
  \includegraphics[width=.9\linewidth]{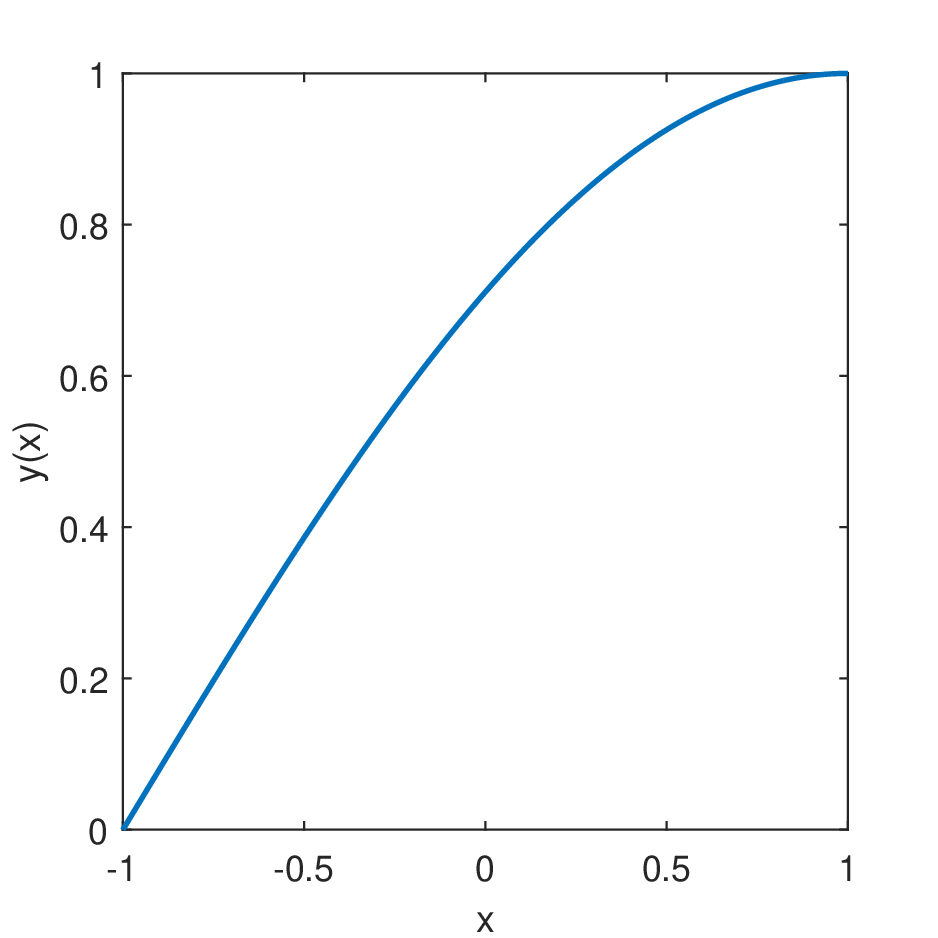}
 \caption{ }
  \label{fig:2DVsub1}
\end{subfigure}%
\begin{subfigure}{.5\textwidth}
  \centering
  \includegraphics[width=.9\linewidth]{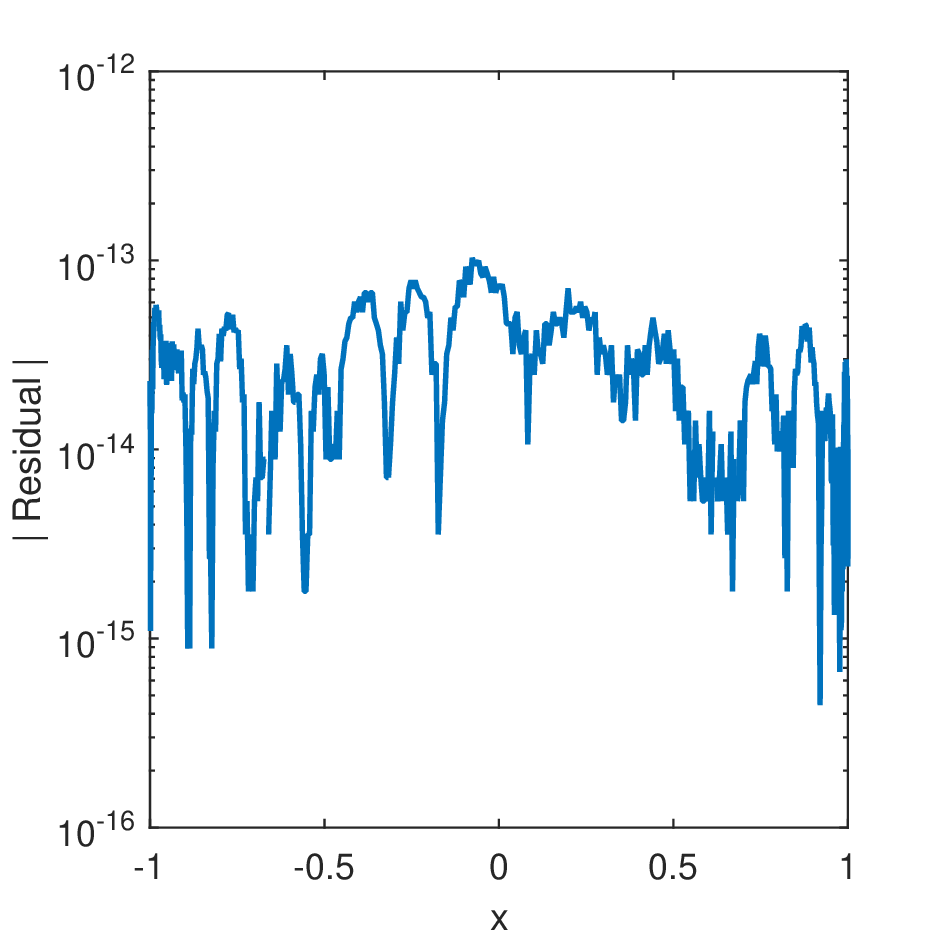}
 \caption{ }
  \label{fig:2DVsub2}
\end{subfigure}
\caption{Solution and error of (\ref{eq:TwoDViscous}) calculated at $\alpha = 1$ and $R_{e} = 10$ using $\mathcal{L}_{4}$ from (\ref{eq:aux_linear_operators})  discretised over $2^9$ grid points using the GHAM.}
\label{fig:2DVCaption}
\end{figure}
\begin{figure}
\centering
  \includegraphics[width=.55\linewidth]{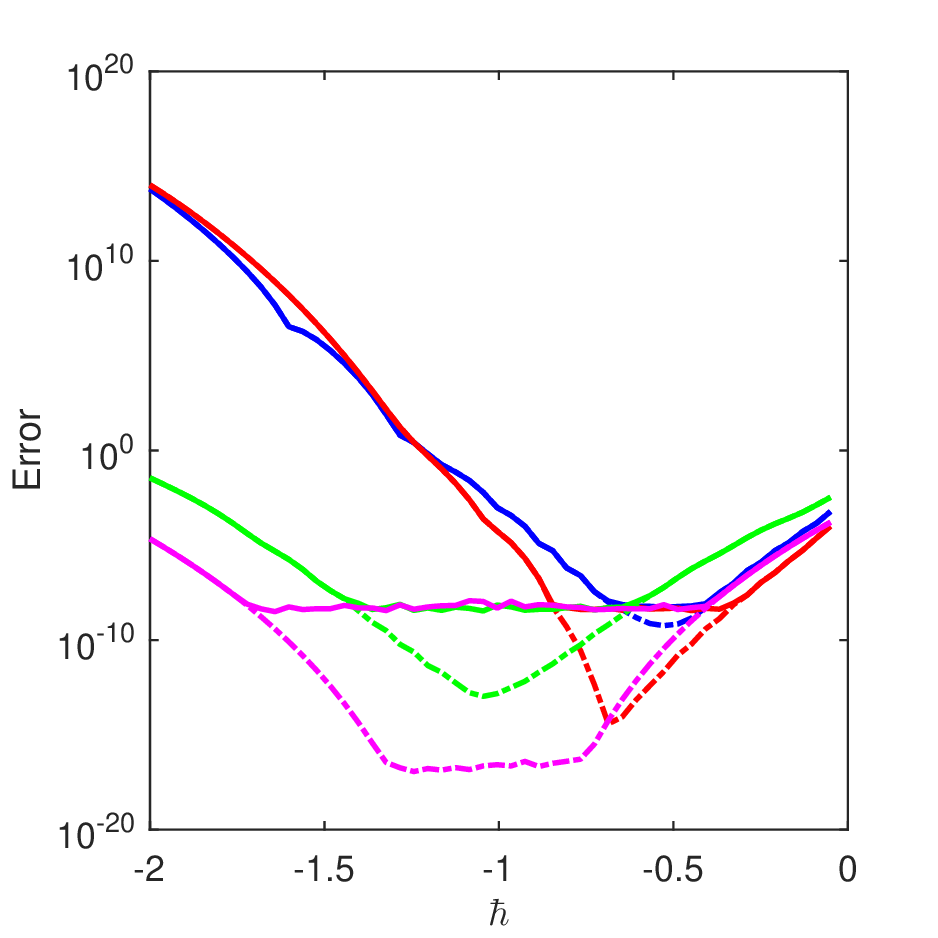}
\caption{Error against $\hbar$ for solutions of (\ref{eq:TwoDViscous}), as calculated using the SHAM (solid lines) and the GHAM (dotted lines). The blue lines correspond to $\mathcal{L}_{1}$, the red to $\mathcal{L}_{2}$, the green to $\mathcal{L}_{3}$ and the magenta line corresponds to $\mathcal{L}_4$, with all operations truncated after $25$ iterations.}
\label{fig:TwoDViscousHBar}
\end{figure}
All four choices of linear operator can be used to construct solutions to the motivating nonlinear problem, as shown in Figure \ref{fig:2DVCaption}. Each operator will exhibit different convergence properties, with each scheme differing in the complexity of constructing solutions numerically and secondly, the operators' resulting convergence properties. Figure \ref{fig:TwoDViscousHBar} explores the relationship between the choice of $\hbar$ and the error that results for each of the linear operators, where the error is defined as the integral over the residual, i.e. $\int_{0}^{1} |R(x)|_{2} dx$, where $|\cdot|_{2}$ denotes the $L_{2}$ norm. The optimal value of $\hbar$ can be considered to be the value of $\hbar$ where the error is minimised, for a given number of iterations. However, in this example we see that as the number of iterations increases, the  error hits a threshold value, below which the solution from the SHAM will not improve. For solutions of (\ref{eq:TwoDViscous}), both the SHAM and the GHAM show that the schemes share the same convergence properties for each linear operator, with the only point of difference being that the implementation of the SHAM used to calculate these results is unable to calculate solutions with greater accuracy than $\mathcal{O}(10^{-10})$. As might be expected, as the operator $\mathcal{L}$ is modified to take a closer form to the full nonlinear equation, the observed convergent region at $25$ iterations widens, and, generally, yields higher accuracy solutions. The exception to this is that while $\mathcal{L}_{2}$ has a smaller convergent region, at its optimal $\hbar$---approximately $0.75$---the calculated error is markedly lower than that calculated using $\mathcal{L}_{3}$ at its optimal $\hbar$. Based upon results from other nonlinear equations, this pattern generally holds---that as long as the auxiliary linear operator $\mathcal{L}$ reflects, in some basic sense, the original dynamics of the nonlinear equation being solved, then the scheme will be convergent. However by more closely replicating the full nonlinear equation greater convergence properties can be observed---which will have implications for the computational cost of constructing a solution to the nonlinear problem at hand. 

A notable property of the error, as exhibited within Figure \ref{fig:TwoDViscousHBar}, is that it is broadly convex with respect to $\hbar$, a result that has been confirmed through other testing. This convexity is advantageous, as it significantly simplifies the process of searching for the $\hbar$ at which the solutions converges the fastest, $\hbar_{\text{opt}}$. This search can be performed at low iterations, with the solution at $\hbar_{\text{opt}}$ then being iterated upon until a desired error tolerance is satisfied.

\begin{figure}
\centering
\begin{subfigure}{.5\textwidth}
  \centering
  \includegraphics[width=.9\linewidth]{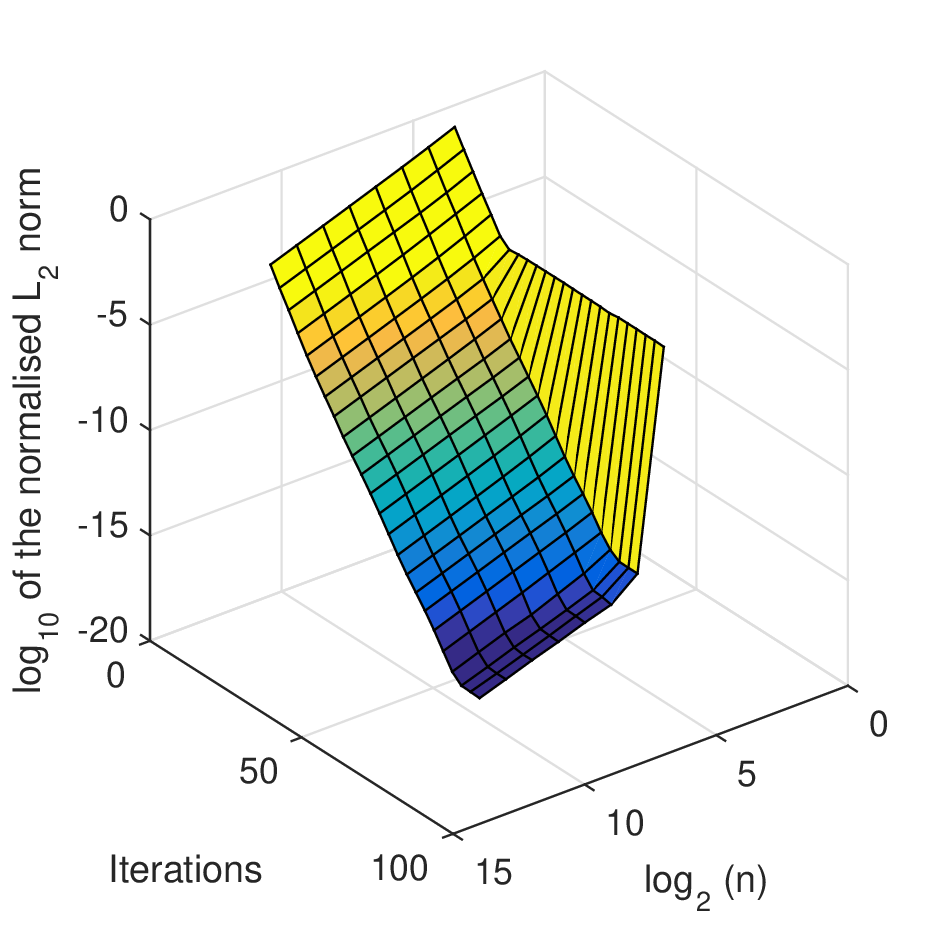}
  \caption{$\mathcal{L}_{1}$}
  \label{fig:LaminarTransVsub1}
\end{subfigure}%
\begin{subfigure}{.5\textwidth}
  \centering
  \includegraphics[width=.9\linewidth]{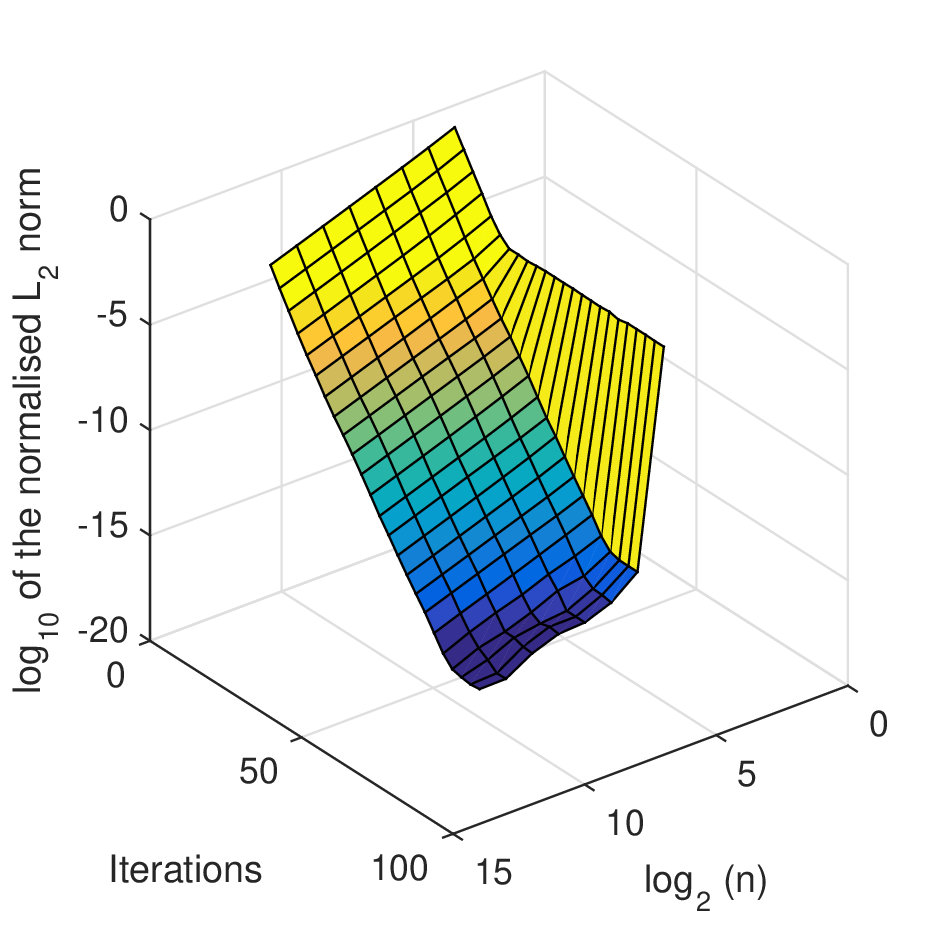}
  \caption{$\mathcal{L}_{2}$}
  \label{fig:LaminarTransVsub2}
\end{subfigure}
\begin{subfigure}{.5\textwidth}
  \centering
  \includegraphics[width=.9\linewidth]{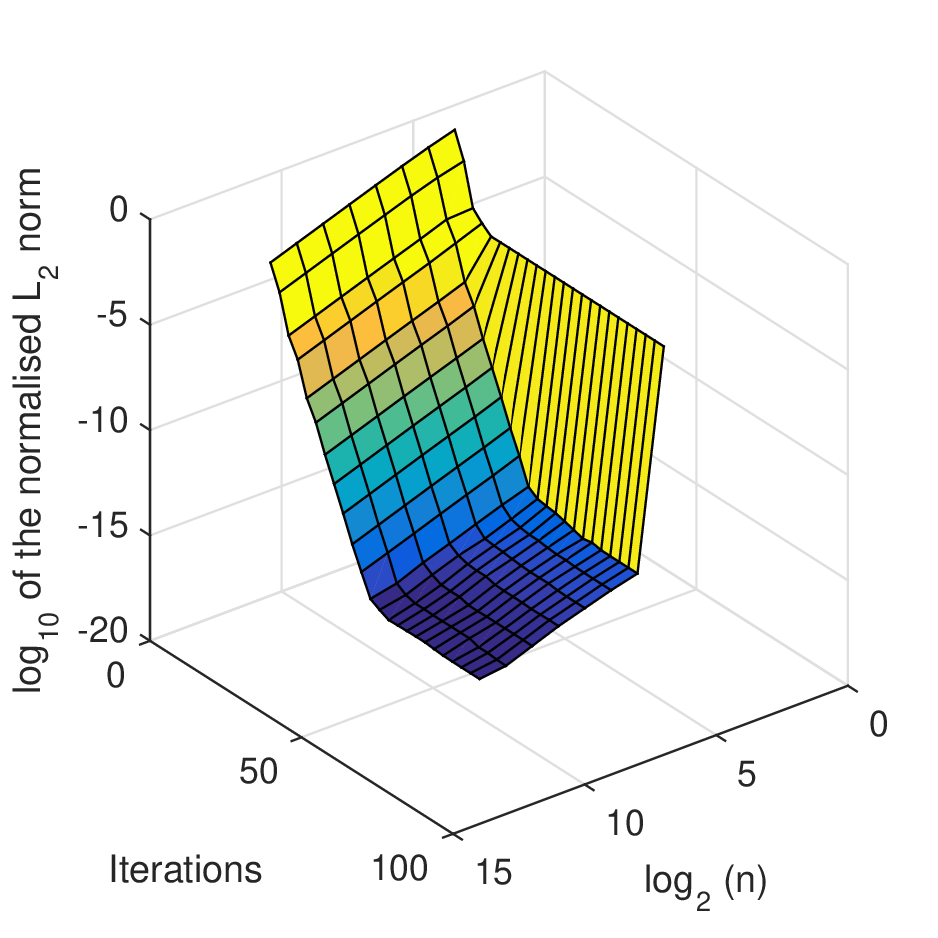}
  \caption{$\mathcal{L}_{3}$}
  \label{fig:LaminarTransVsub3}
\end{subfigure}%
\begin{subfigure}{.5\textwidth}
  \centering
  \includegraphics[width=.9\linewidth]{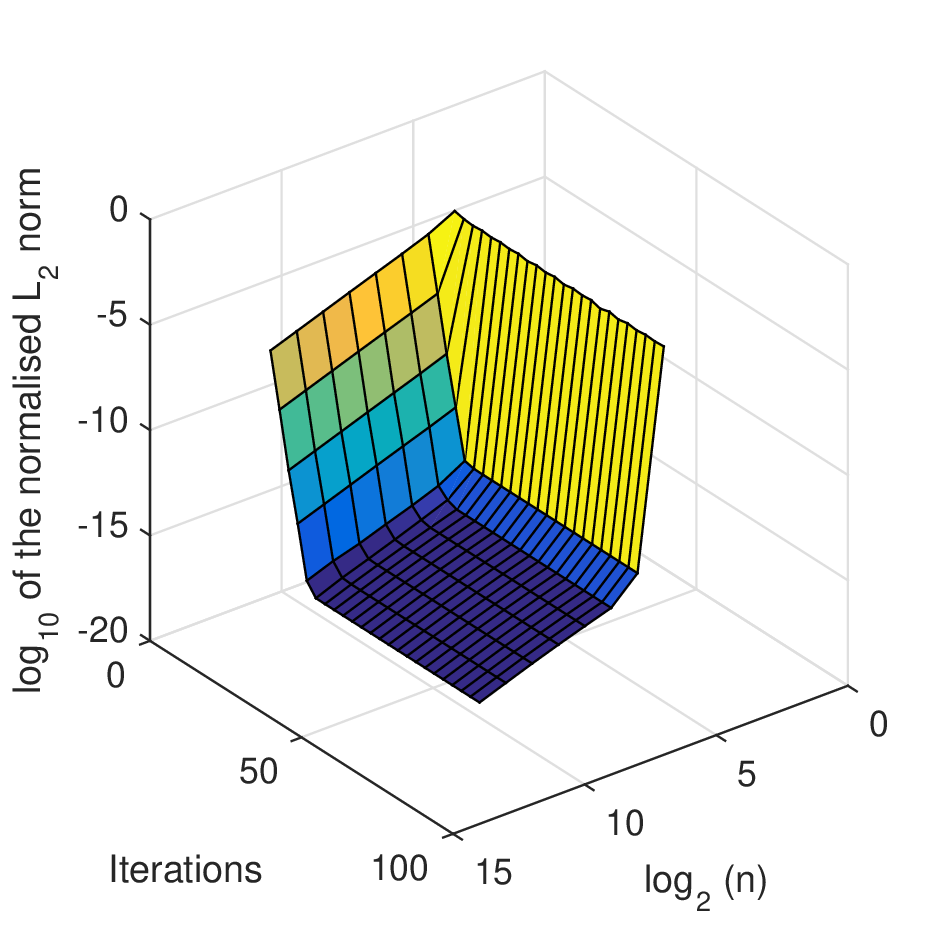}
  \caption{$\mathcal{L}_{4}$}
  \label{fig:LaminarTransVsub4}
\end{subfigure}
\caption{Calculated error at $\hbar_{\text{opt}}$ for (\ref{eq:TwoDViscous}) using the GHAM, as a function of the number of iterations and spatial resolution $n$. For these calculations $H_{a} = 1$ and $R_{e} = 10$.}
\label{fig:LaminarTransVCaption}
\end{figure}

\begin{figure}
\centering
\begin{subfigure}{.5\textwidth}
  \centering
  \includegraphics[width=.9\linewidth]{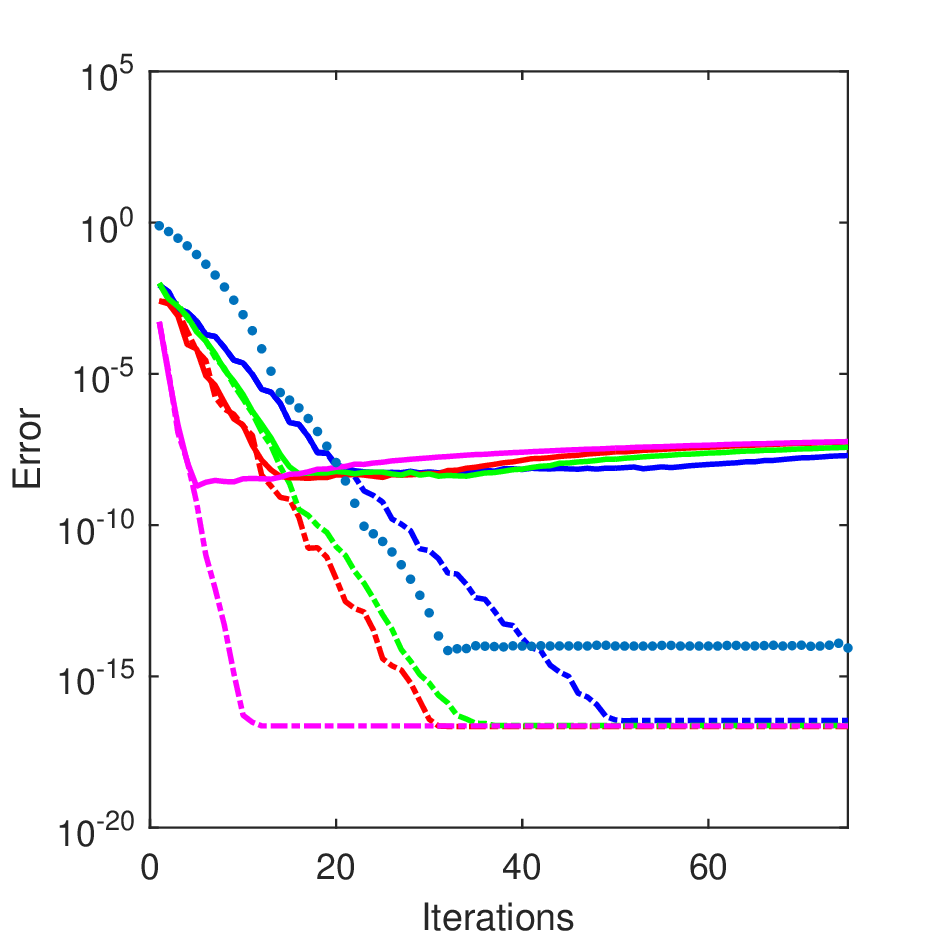}
  \caption{$2^8$ points}
  \label{fig:LaminarTransVsub1SHAMvGHAM}
\end{subfigure}%
\begin{subfigure}{.5\textwidth}
  \centering
  \includegraphics[width=.9\linewidth]{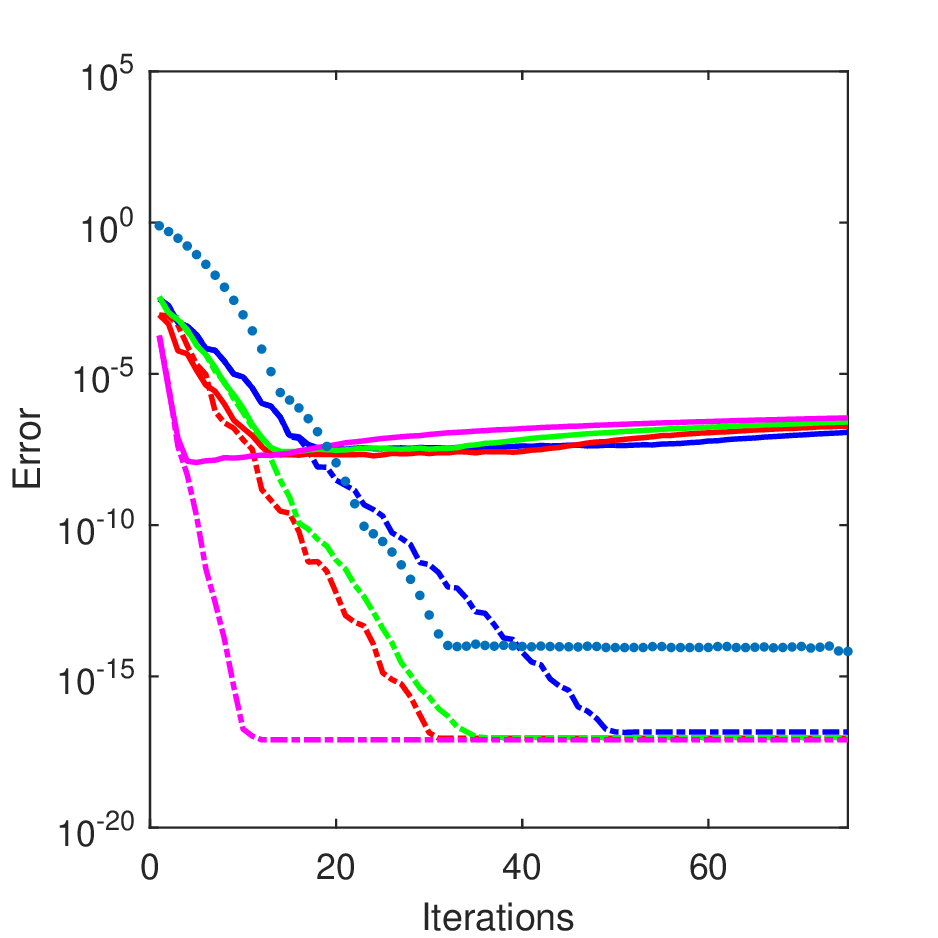}
  \caption{$2^9$ points }
  \label{fig:LaminarTransVsub2SHAMvGHAM}
\end{subfigure}
\caption{Impact of the number of iterative steps upon the numerical error when solving (\ref{eq:TwoDViscous}) for varying choices of the auxiliary linear operator $\mathcal{L}$. Solid lines correspond to the SHAM, dashed denote the GHAM, and the blue dots are Newton Iteration upon a Gegenbauer discretisation. Of the solid lines, blue represents $\mathcal{L}_1$, Red is $\mathcal{L}_2$, Green is $\mathcal{L}_3$ and Magenta is $\mathcal{L}_4$. All solutions were calculated at the optimal $\hbar$ for each choice of $\mathcal{L}$.}
\label{fig:LaminarTransVCaptionSHAMvGHAM-1}
\end{figure}

Figure \ref{fig:LaminarTransVCaption} presents an analysis of the error, as driven by both the convergence of the nonlinear problem, and from the linear sub-problems. As the grid resolution increases beyond a threshold value of $n = 2^4$, the error of the solution becomes nominally independent of resolution. This results from the resolution passing the threshold for the spectral convergence of the Gegenbauer method for the linear problem, so that beyond this point changes in the error are driven by the process of solving the nonlinear problem. This can be seen in the change of the convergence of the solution as the number of iterations increases. As long as the threshold resolution has been reached, the primary driver of the convergence of the solution is the number of iterations, up until the point where the numerical solution has converged, within machine precision, to the true solution. While the different auxiliary linear operators of~\autoref{eq:aux_linear_operators} exhibit different rates of convergence, as will be shown further in Equation \ref{fig:LaminarTransVsub2SHAMvGHAM}, the spectral convergence behaviour are still shared across all choices of $\mathcal{L}$.

\begin{figure}
\centering
\begin{subfigure}{.5\textwidth}
  \centering
  \includegraphics[width=.9\linewidth]{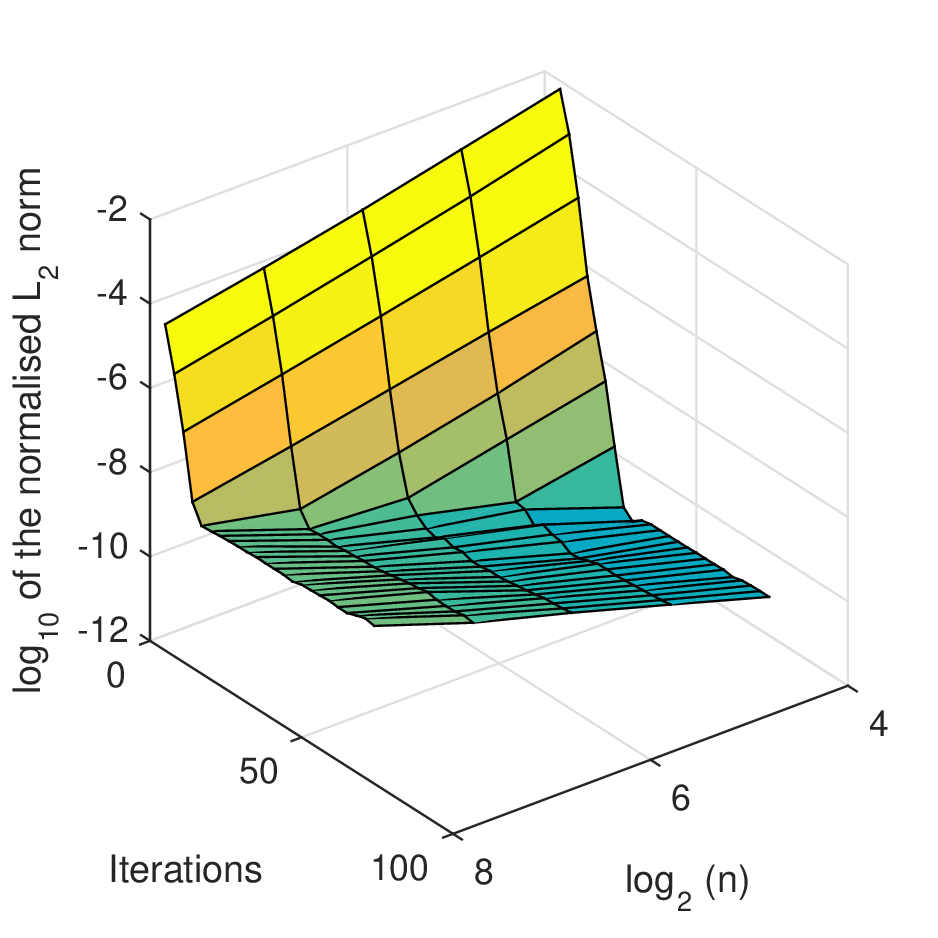}
  \caption{$\mathcal{L}_{1}$}
  \label{fig:LaminarTransVsub1SHAM}
\end{subfigure}%
\begin{subfigure}{.5\textwidth}
  \centering
  \includegraphics[width=.9\linewidth]{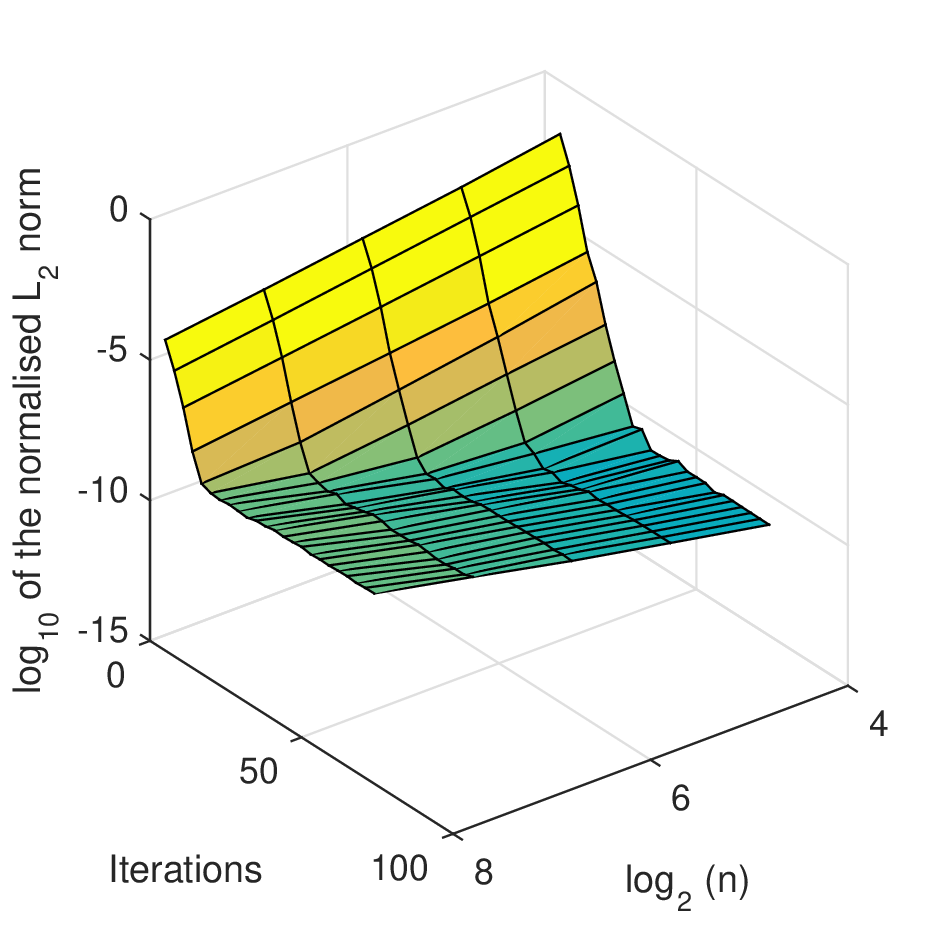}
  \caption{$\mathcal{L}_{2}$}
  \label{fig:LaminarTransVsub2SHAM}
\end{subfigure}
\begin{subfigure}{.5\textwidth}
  \centering
  \includegraphics[width=.9\linewidth]{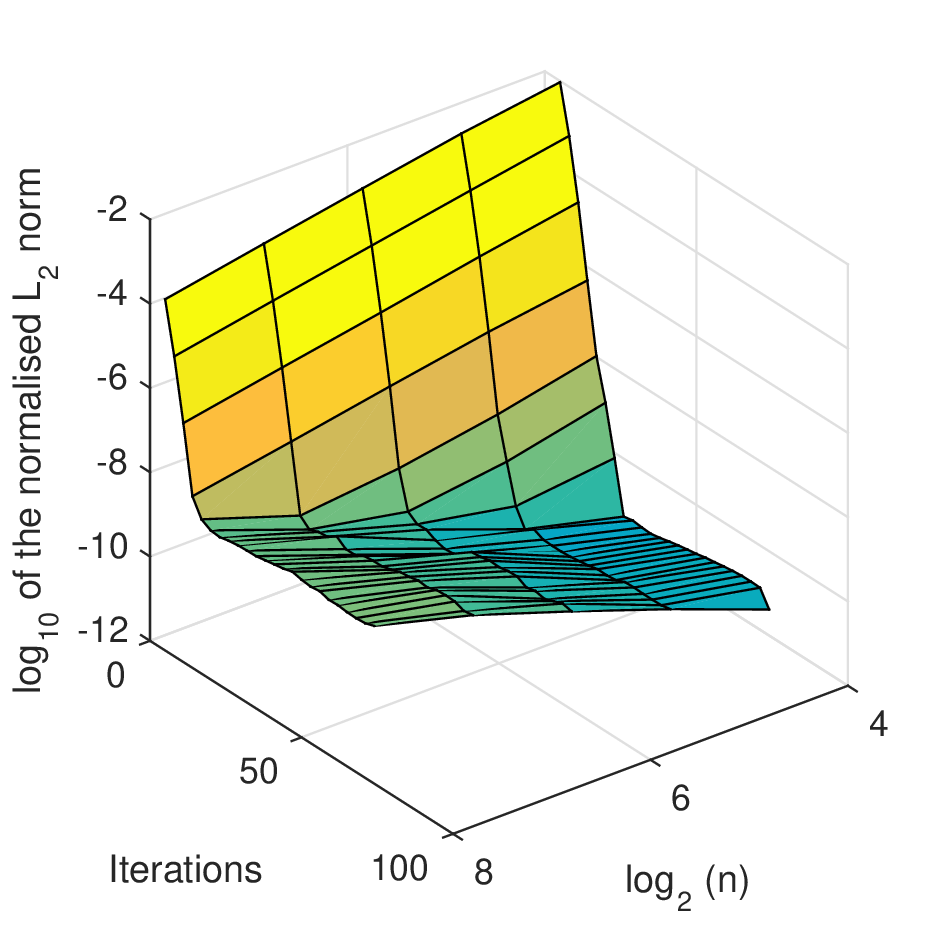}
  \caption{$\mathcal{L}_{3}$}
  \label{fig:LaminarTransVsub3SHAM}
\end{subfigure}%
\begin{subfigure}{.5\textwidth}
  \centering
  \includegraphics[width=.9\linewidth]{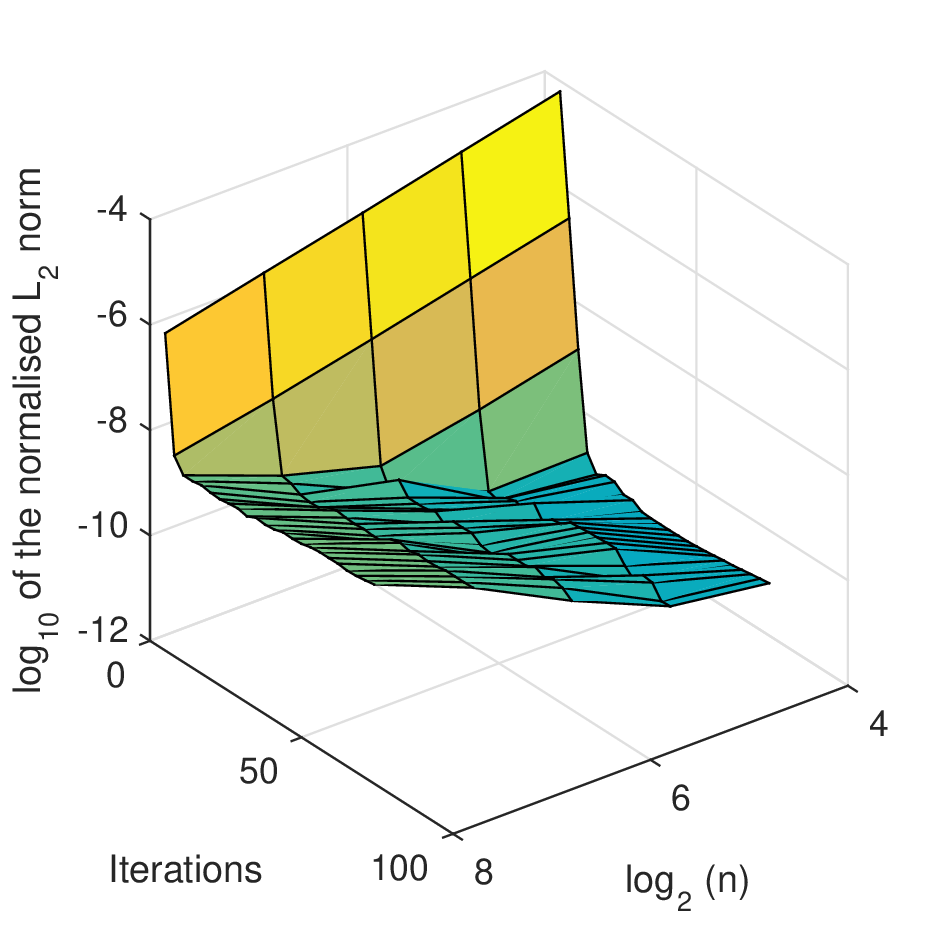}
  \caption{$\mathcal{L}_{4}$}
  \label{fig:LaminarTransVsub5SHAM}
\end{subfigure}
\caption{Calculated error at $\hbar_{\text{opt}}$ for (\ref{eq:TwoDViscous}) using the SHAM, as a function of the number of iterations and spatial resolution $n$. For these calculations $H_{a} = 1$ and $R_{e} = 10$.}
\label{fig:LaminarTransVCaptionSHAM}
\end{figure}

To explore how the dynamics of the GHAM differ from those exhibited by the SHAM, the results from Figure \ref{fig:LaminarTransVCaption} were recreated using the SHAM in Figure \ref{fig:LaminarTransVCaptionSHAM}, where the linear matrix equations are discretised using Chebyshev-collocation matrices. While the convergence results of the SHAM and the GHAM largely seem comparable, there are some fundamental differences between the performances of the two techniques. Given sufficient iterations, the results from the SHAM converge for any spatial resolution, as compared to the GHAM which requires a threshold number of grid points in order to resolve a solution. As the grid resolution increases the accuracy of the GHAM increases, while, somewhat counter intuitively, the accuracy of the SHAM actually decreases. This can be attributed to the Chebyshev collocation matrices that underpin the SHAM becoming singular as the grid resolution increases, which in turn increases the inherent errors in attempting to solve these systems numerically. The factors underpinning this phenomena have been discussed by various authors \citep{Bayliss1995, Belfert1997, Breuer1992, Peyret2002, Rothman1991, Trefethen1989}, with the main cause being the accumulation of roundoff errors in the calculation of the differentiation matrices, and the effect of the decreased spacing between grid points, especially as $|x| \to 1$. This also explains why the SHAM is unable to match the accuracy exhibited by the GHAM---the latter of which converges to machine precision. 

To this point, the convergence of the GHAM and the SHAM has been considered in the context of their accuracy, and the convergence of the residuals of the calculated solutions as a function of the grid resolution and the number of iterations employed. However, the primary advantage of the sparse matrix formulation of the GHAM is its low theoretical scaling of computational time with respect to the grid resolution. Figure \ref{fig:LaminarTransVCaptionSHAMvGHAM-1} considers the performance of both the GHAM and the SHAM at $\hbar_{\text{opt}}$ for solving (\ref{eq:TwoDViscous}), subject to a selection of linear operators.

For the GHAM solver, the solution in terms of $\mathcal{L}_4$ significantly outperforms all other schemes, converging to machine precision in only $20$ iterations. While its pre--eminence is unsurprising, given how closely $\mathcal{L}_4$ approximates the full nonlinear problem, the degree to which it outperforms the other choices of $\mathcal{L}$ is surprising. In contrast to this, the remaining choices for $\mathcal{L}$ all take approximately three times as many iterations in order to converge---however all three still clearly exhibit spectral convergence. Interestingly, $\mathcal{L}_{2}$ slightly outperforms $\mathcal{L}_{1}$ and $\mathcal{L}_{3}$, even though $\mathcal{L}_{3}$ can be considered to be a more accurate linear representation of the nonlinear problem. These differences are maintained when considering the computational time required for all the solutions.

The relative performance of each of the linear operators was maintained within the SHAM. However, there is a marked difference in the minimum error exhibited by the GHAM, as compared to the SHAM. While the GHAM is able to converge to machine precision, the SHAM cannot resolve solutions beyond a threshold error, and that threshold is not constant as the resolution is changed.

\begin{figure}
\centering
\begin{subfigure}{.5\textwidth}
  \centering
  \includegraphics[width=.9\linewidth]{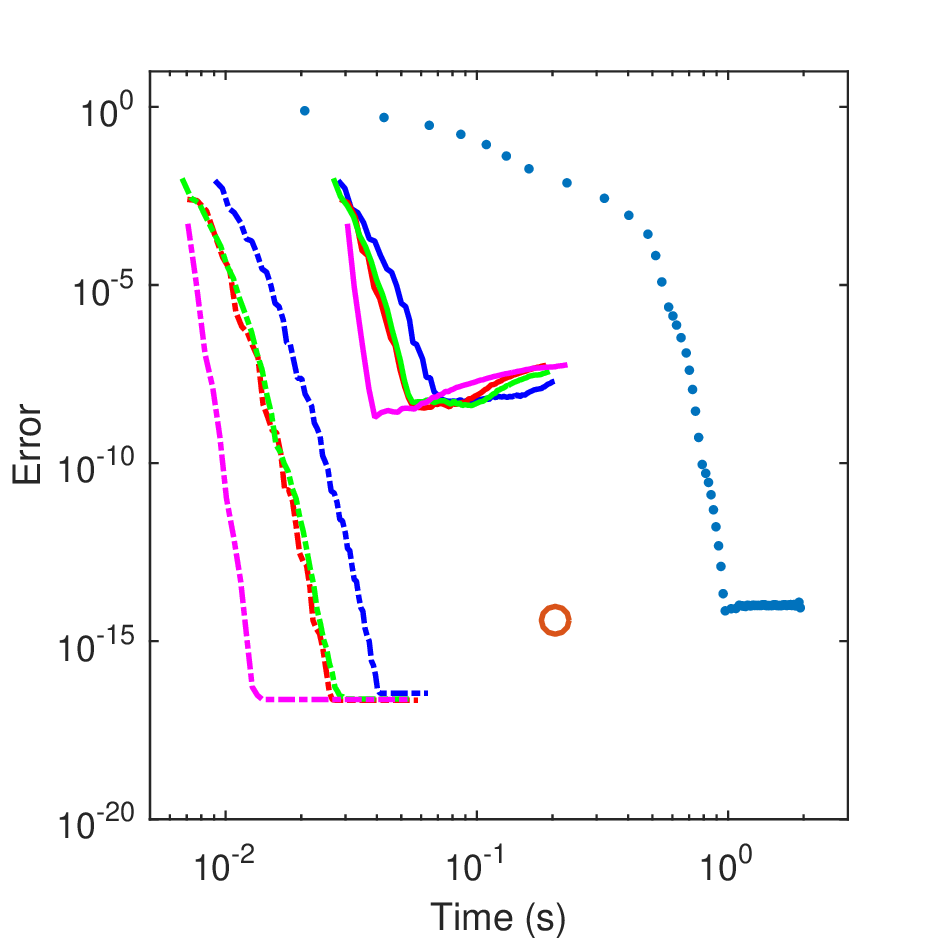}
  \caption{$2^8$ points}
  \label{fig:LaminarTransVsub1SHAMvGHAMvNewton}
\end{subfigure}%
\begin{subfigure}{.5\textwidth}
  \centering
  \includegraphics[width=.9\linewidth]{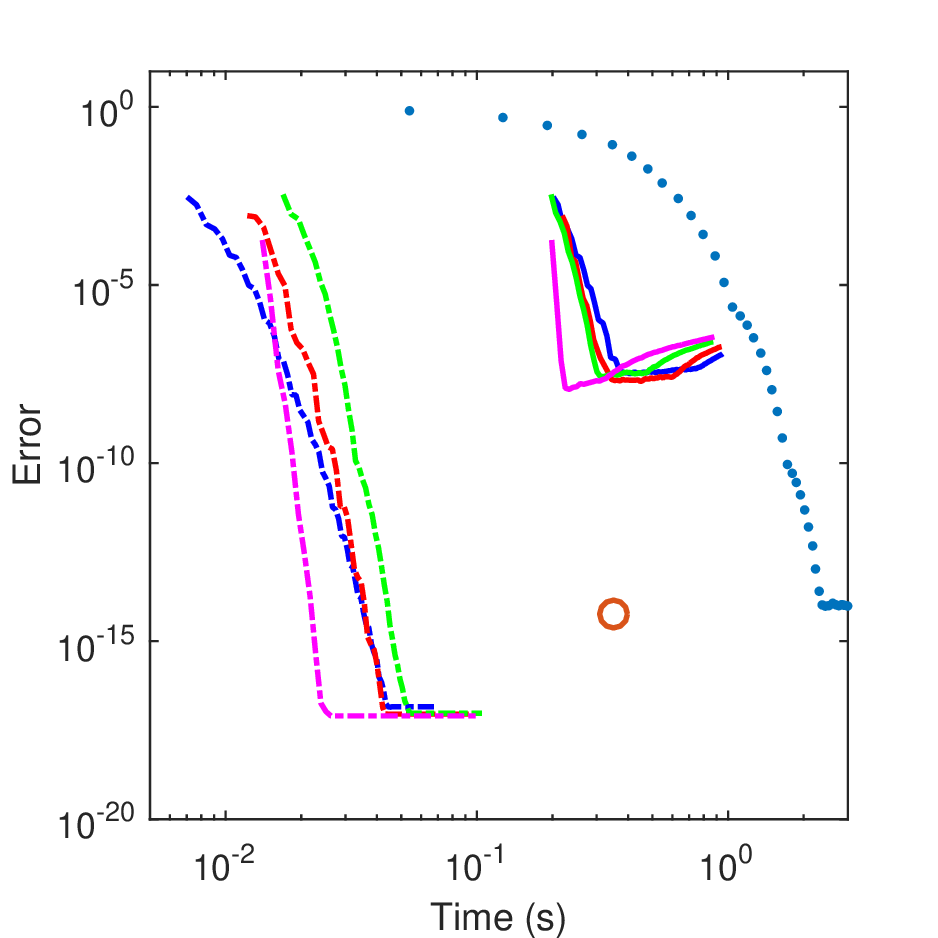}
  \caption{$2^9$ points}
  \label{fig:LaminarTransVsub2SHAMvGHAMvNewton}
\end{subfigure}
\caption{Computational cost and the error of solving (\ref{eq:TwoDViscous}) for varying numbers of iterations, following Figure \ref{fig:LaminarTransVCaptionSHAMvGHAM-1} with the inclusion of the computational cost and error using MATLAB's BVP4c routine, as represented by the red circle.}
\label{fig:LaminarTransVCaptionSHAMvGHAMvNewton}
\end{figure}

To explore the relative numerical performance of these schemes, two other numerical approaches were taken to solve the motivating equation. The first was Newton Iteration based upon a Gegenbauer discretisation (as described by \cite{Olver2013}), in order to separate our the contribution to the numerical efficiency from the nonlinear approach, and then the linear solver that both approach is built upon. The second comparison was to MATLAB's `BVP4C' routine. The latter comparison should theoretically be unfavourable towards GHAM, as BVP4C is a highly optimised routine written in compiled C, as compared to the GHAM codebase which is in terms of uncompiled MATLAB code. However, with both approaches, as shown in Figure \ref{fig:LaminarTransVCaptionSHAMvGHAMvNewton}, GHAM for the optimal choice of $\mathcal{L}$ and $\hbar$ is multiple orders of magnitude faster than the comparison approaches. Even the dense matrices of SHAM outperform Newton iteration across all tested resolutions, due to Newton iteration necessitating the construction of a new dense matrix operator at each step.

As the resolution of the grid discretisation is doubled, the time required to solve using both the SHAM and Newton--Iteration increases by almost an order of magnitude, whereas there is only an approximate doubling in the computational time required when the calculations are performed with the GHAM. This result has been replicated across other resolutions, and stems from the rate of growth of the number of elements in the sparse matrix operators used within the GHAM, as compared to all the other tested techniques. 

\begin{figure}
\centering
\begin{subfigure}{.5\textwidth}
  \centering
  \includegraphics[width=.8\linewidth]{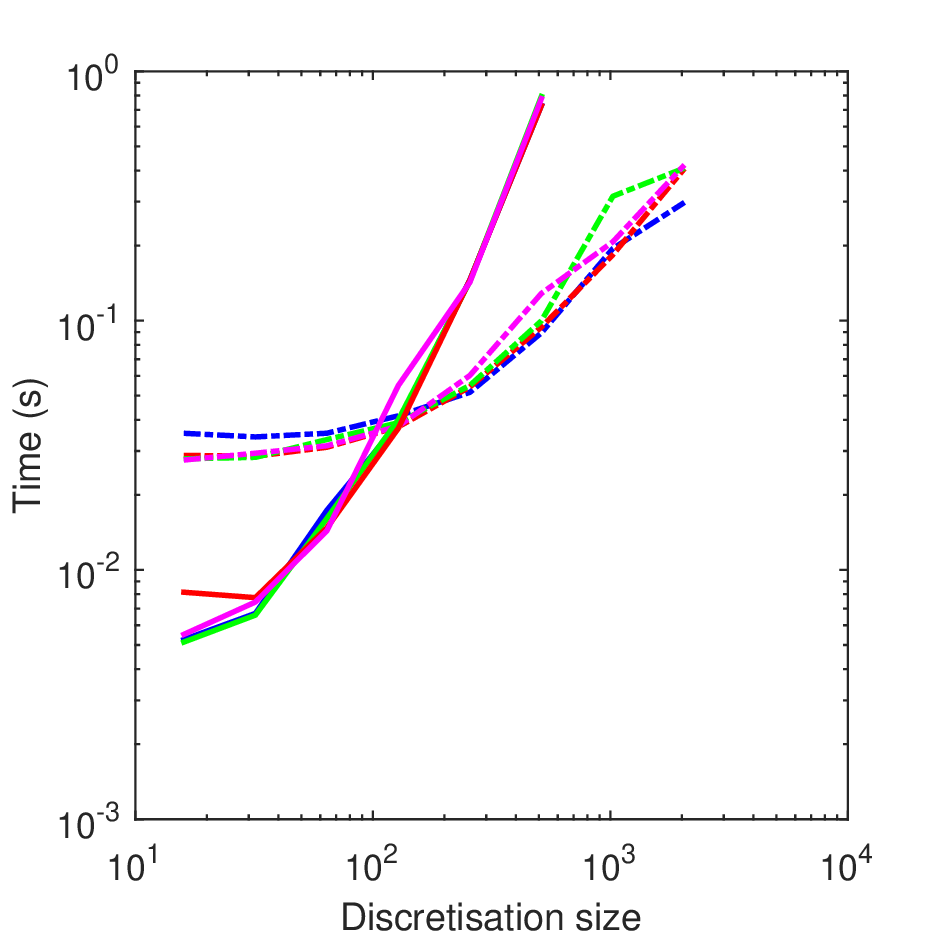}
  \caption{Grid resolution and computational time.}
  \label{fig:TwoDViscoussub1SHAMvGHAM}
\end{subfigure}%
\begin{subfigure}{.5\textwidth}
  \centering
  \includegraphics[width=.8\linewidth]{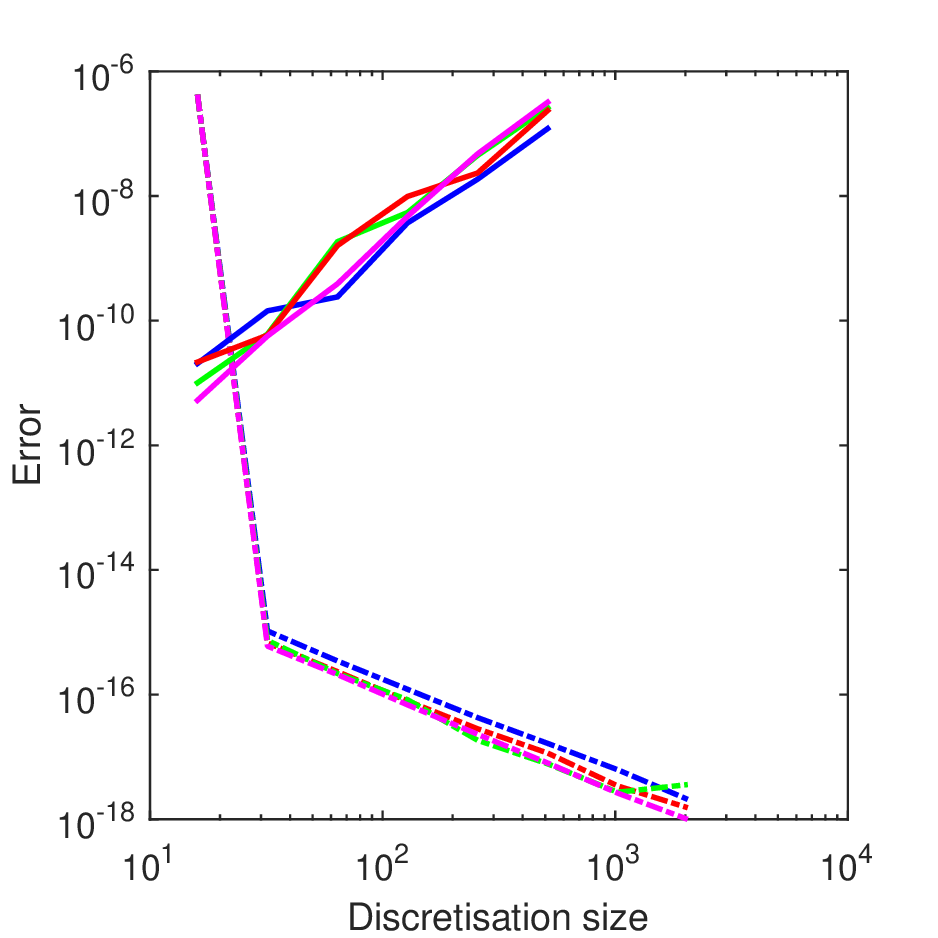}
  \caption{Grid resolution and error.}
  \label{fig:TwoDViscoussub2SHAMvGHAM}
\end{subfigure}
\caption{Scaling of computational time for solutions of (\ref{eq:TwoDViscous}) calculated by taking $75$ iterations of the SHAM (solid lines) and the GHAM (dotted lines). Results for the GHAM are presented for $n \in [2^6,2^{14}]$, whereas the SHAM was only calculated over $n \in [2^6,2^9]$, as the Chebyshev collocation matrices within the SHAM become singular after this point.}
\label{fig:TwoDViscoussVCaptionSHAMvGHAM}
\end{figure}

The relative differences between the numerical performance of these schemes is maintained across varying resolutions, with the GHAM significantly outperforming the SHAM at nearly all resolutions tested within Figure \ref{fig:TwoDViscoussVCaptionSHAMvGHAM}. For small resolutions, the cost of establishing the matrix system dominates, which results in the SHAM exhibiting a relatively lower computational time. However, as $n$ increases, the cost of solving the system begins to dominate, leading to the GHAM starting to significantly outperforming the SHAM. An additional consideration of note is that unlike the GHAM, the SHAM could not be considered for grid discretisations larger than $2^9$ points, as the matrices beyond this point become singular. More broadly, as $n$ approaches $2^9$, the determinant of the linear matrix in the SHAM approaches zero, which likely explains the growth in the error shown within Figure \ref{fig:TwoDViscoussub2SHAMvGHAM}.

To delve into the factors that drive the computational cost for the GHAM, the process of solving (\ref{eq:TwoDViscous}) can be separated out into its constituent parts. The process first involves discretising the system using Gegenbauer polynomials, followed by solving the matrix system at each iteration, converting between Chebyshev space and real space and evaluating the derivatives $\{y^{(1)},y^{(2)},y^{(3)},y^{(4)} \}$. 

To verify the previously outlined theoretical results for how the GHAM scales with the number of iterations, a regression to the form of $t = CI^{S}$ was conducted, where $I$ is the number of iterations and $C$ and $S$ are proportionality constants. This regression was based upon the computational time to calculate a solution between $40$ and $80$ iterations, in order to isolate out the effect of the iterative process from the first step of the scheme, and to examine how the cost of constructing the inhomogeneous component scales, as it should become more costly as the number of iterations increases. For each linear operator $C$ scales with $n$, which reflects the cost of solving the matrix operations after $LU$ decomposition, and how that grows with the number of grid points in the discretisation. However, $S$ is approximately $1$ for all linear operators, across all grid resolutions, reflecting that the cost of the scheme is entirely driven by the $LU$ matrix operations, rather than any other components of the iterative process.

\begin{figure}
\centering
\begin{subfigure}{.5\textwidth}
  \centering
  \includegraphics[width=.9\linewidth]{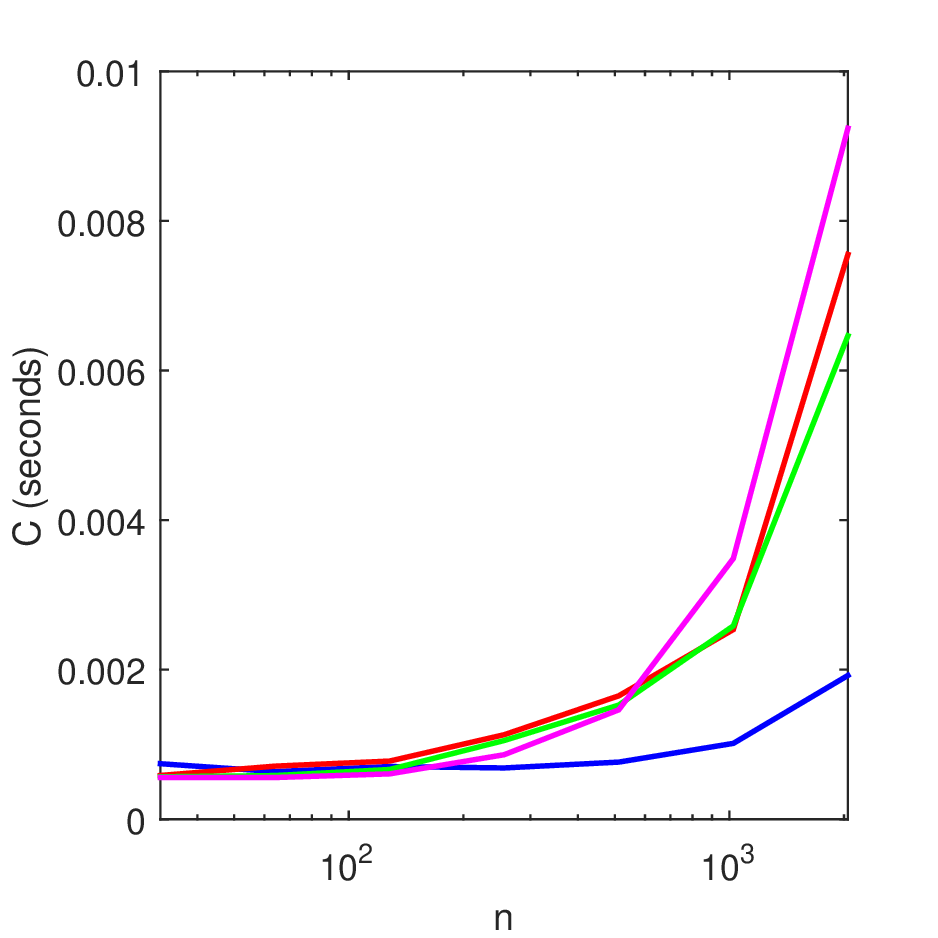}
  \caption{$\mathcal{L}_{1}$}
  \label{fig:C_scaling}
\end{subfigure}%
\begin{subfigure}{.5\textwidth}
  \centering
  \includegraphics[width=.9\linewidth]{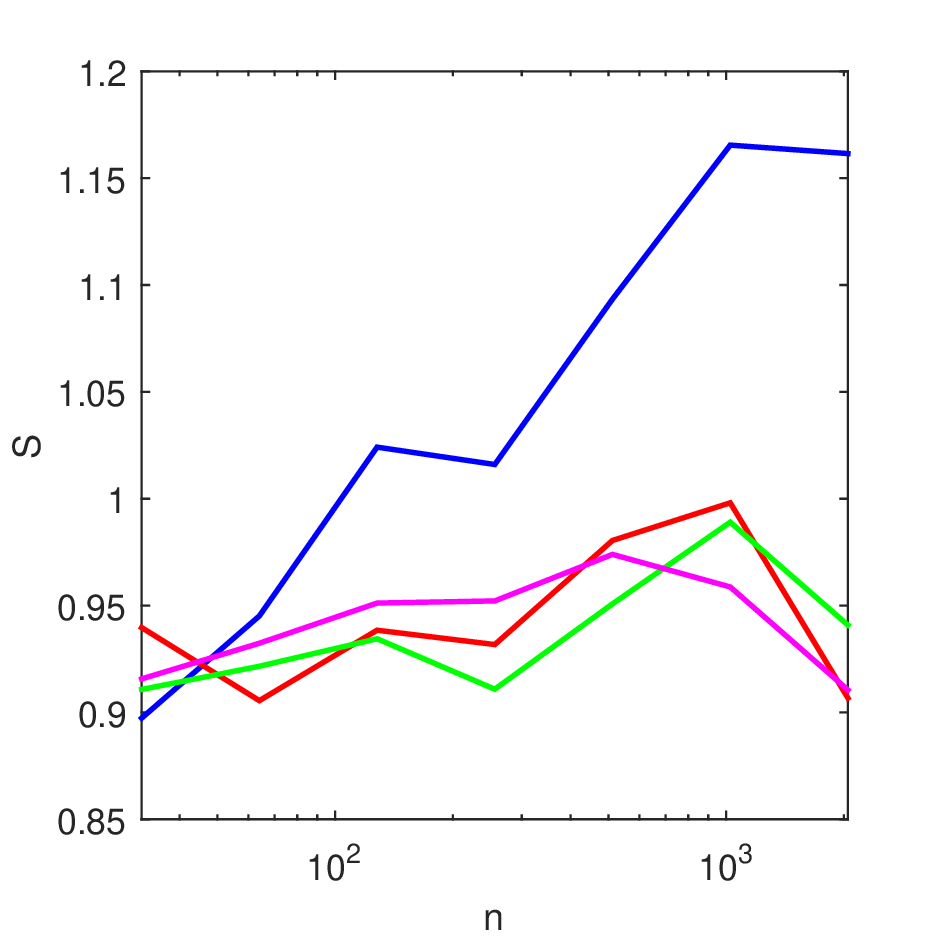}
  \caption{$\mathcal{L}_{2}$}
  \label{fig:S_scaling}
\end{subfigure}
\caption{Scaling coefficients for the GHAM, assuming that the computational cost scales as $T = C I^{S}$. One again, blue, red, green and magenta represent $L_{1}$, $L_{2}$, $L_{3}$ and $L_{4}$ respectively.}
\label{fig:scaling}
\end{figure}

\begin{figure}
\centering
  \includegraphics[width=.65\linewidth]{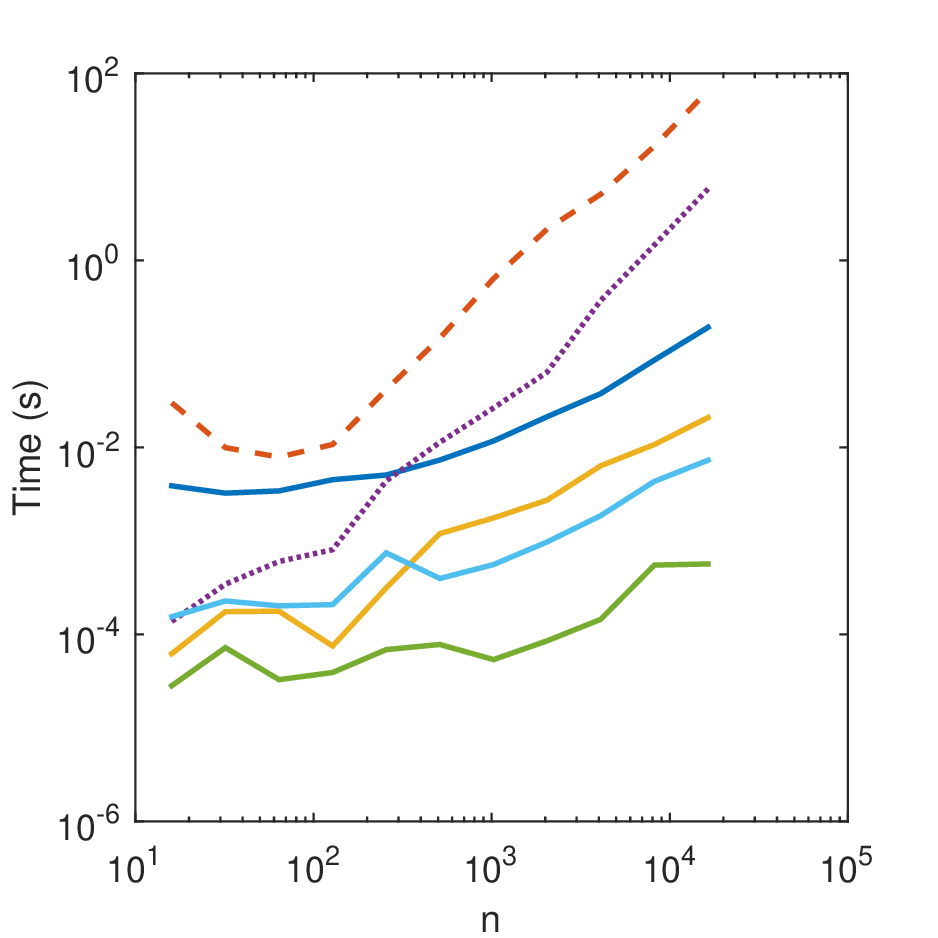}
  \caption{Scaling of computational cost of constructing solutions for~\autoref{eq:TwoDViscous} with $n$. Blue: set up cost for establishing the matrix problem in the context of the GHAM. Red, dashed: LUPQR decomposition, which only needs to occur once at the beginning of the iterative process. Yellow: solving the matrix system using the LUPQR decomposition. Green: transform between real and Chebyshev space. Light blue: calculating derivatives up to fourth--order. Purple, dotted: solving the matrix system using a direct matrix inverse, which is not used within the algorithm.} 
  \label{fig:TwoDViscousSolveCost}
\end{figure}

Clearly the dominant component of the computational cost stems from the grid resolution, rather than the number of iterations required to solve the problem. To further explore the scaling of computational cost, the results for the GHAM in Figure \ref{fig:TwoDViscoussub1SHAMvGHAM} can be decomposed into its constituent parts, the results for which are presented in Figure \ref{fig:TwoDViscousSolveCost}. The dominant contributors to solving this numerical problem are the cost of performing the $LUPQR$ decomposition; solving the matrix system after the $LU$-decomposition has been employed; evaluating the first four spatial derivatives; and the transform between real and Chebyshev space and visa-versa. While the cost of performing the $LUPQR$ decomposition would appears to be unfavourably high with respect to the direct cost of solving the matrix systems, as shown in Figure \ref{fig:TwoDViscousSolveCost}, it must be stressed that directly solving the matrix would be required at every step of the iterative process, whereas the $LUPQR$ decomposition only has to be performed once, requiring $\mathcal{O}(n^{1.305})$ time at the beginning of the iterative process. 
The other components of Figure \ref{fig:TwoDViscousSolveCost} do need to be performed at each step of the iterative process, and it is their cost that influences how the scheme scales with the number of iterations---behaving in aggregate as an $\mathcal{O}(n^{1.05})$ process. This quasi--linear scaling is entirely a product of the low fill--in density of the matrix operators that make up the GHAM.

Based upon the results from (\ref{eq:TwoDViscous}) it can be said that---at least for the tested fourth--order systems---the scheme scales with $\mathcal{O}\left( n^{1.305} + m n^{1.05} \right)$, where $n$ is the spatial resolution and $m$ is the total number of iterative steps required to both find $\hbar_{\text{opt}}$, and then to refine the solution to a prescribed error tolerance. This compares favourably to the SHAM, as solving dense matrix systems is limited by the theoretical scaling $\mathcal{O}(n^{2.373})$. Interestingly, while it was predicted in the preceeding section that there would be nonlinear scaling with the number of iterations, stemming from the $\mathcal{O}( m^2 n)$ multiplication operations, in practice this term is dominated by the costs involved in solving the matrix system, which scale linearly with iterations. This result has been tested up to $150$ iterations, at which point the system still appears to scale linearly with $m$.

\section{Conclusion}

We have presented a new, spectrally accurate scheme for solving steady, nonlinear boundary value problems, that uses the Homotopy Analysis Method to convert the original problem into an infinite sequence of linear differential equations, which are in turn solved upon a Gegenbauer discretisation. Unlike other numerical approaches for steady, nonlinear boundary value problems, the GHAM exhibits quasi--linear growth in time with respect to the grid resolution, and the number of iterations, due to the iterative process being constructed in terms of a single inversion of a sparse, almost--banded matrix. In the context of the examined test problem, we showed that this scaling does not impact upon the accuracy exhibited by the scheme, and that it in fact is able to resolve more accurate solutions than the comparable SHAM.

While the scheme has proven accurate for one--dimensional problems, further work is required in order to generalise it to higher--dimensional problems, in order for it to have broader applicability to real world problems. We also believe that there is still scope for improving the computational performance of the scheme. Furthermore, to this point we have not addressed the relevance of the GHAM to stiff differential equations. However, the performance of the scheme, in terms of both its accuracy and its computational cost is highly promising, and as such warrants further exploration.





\bibliographystyle{cbe}
\bibliography{thesis}

\begin{thebibliography}{}

\bibitem[\protect\astroncite{Bayliss et~al.}{1995}]{Bayliss1995}
{\sc Bayliss, A., Class, A., and Matkowsky, B.~J.} 1995.
\newblock Roundoff {E}rror in {C}omputing {D}erivatives {U}sing the {C}hebyshev
  {D}iffereniation {M}atrix.
\newblock {\em Engineering Sciences and Applied Mathematics} 116:380--38.

\bibitem[\protect\astroncite{Belfert}{1997}]{Belfert1997}
{\sc Belfert, B.~D.} 1997.
\newblock Generation of {P}suedospectral {D}iffereniation {M}atrices {I}.
\newblock {\em SIAM Journal on Numerical Analysis} 34:1640--1657.

\bibitem[\protect\astroncite{Breuer and Everson}{1992}]{Breuer1992}
{\sc Breuer, K.~S. and Everson, R.~M.} 1992.
\newblock On the {E}rrors {I}ncurred {C}alculating {D}erivatives {U}sing
  {C}hebyshev {P}olynomials.
\newblock {\em Journal of Computational Physics} .

\bibitem[\protect\astroncite{Chan et~al.}{2017}]{Chan2017}
{\sc Chan, J., Hewett, R., and Warburton, T.} 2017.
\newblock Weight-{A}djusted {D}iscontinuous {G}alerkin {M}ethods: {C}urvilinear
  {M}eshes.
\newblock {\em SIAM Journal on Scientific Computing} 39(6):2395--2421.

\bibitem[\protect\astroncite{Davie and Stothers}{2013}]{Davie2013}
{\sc Davie, A.~M. and Stothers, A.~J.} 2013.
\newblock {I}mproved bound for complexity of {M}atrix {M}ultiplication.
\newblock {\em {P}roceedings of the {R}oyal {S}ociety of {E}dinburgh}
  143A:351--370.

\bibitem[\protect\astroncite{Davis}{2006}]{Davis2006}
{\sc Davis, T.~A.} 2006.
\newblock Direct {M}ethods for {S}parse {L}inear {S}ystems.
\newblock SIAM, Philadelphia.

\bibitem[\protect\astroncite{Davis}{2011}]{Davis2011}
{\sc Davis, T.~A.} 2011.
\newblock Algorithm 915, {S}uisesparseqr: {M}ultifrontal {M}ultithreaded
  {R}ank-{R}evealing {S}parse {QR} {F}actorization.
\newblock {\em ACM Transcations on Mathematical Software} 38:8:1--8:22.

\bibitem[\protect\astroncite{Dhai and Yuan}{1999}]{Dhai1999}
{\sc Dhai, Y.-H. and Yuan, Y.} 1999.
\newblock A {N}onlinear {C}onjugate {G}radient {M}ethod with a {S}trong
  {G}lobal {C}onvergence {P}roperty.
\newblock {\em SIAM Journal of Optimization} 1:177--182.

\bibitem[\protect\astroncite{Gines et~al.}{1998}]{Gines1998}
{\sc Gines, D., Belylkin, G., and Dunn, J.} 1998.
\newblock {LU} {F}actorization of {N}on-standard {F}orms and {D}irect
  {M}ultiresolution {S}olvers.
\newblock {\em Applied and Computational Harmonic Analysis} .

\bibitem[\protect\astroncite{Hayat and Sajid}{2007}]{Hayat2007}
{\sc Hayat, T. and Sajid, M.} 2007.
\newblock {H}omotopy {A}nalysis of {MHD} {B}oundary {L}ayer {F}low of an
  {U}pper-{C}onvected {M}axwell fluid.
\newblock {\em {I}nternational {J}ournal of {E}ngineering {S}cience}
  45:393--401.

\bibitem[\protect\astroncite{He}{2003}]{He2003}
{\sc He, J.-H.} 2003.
\newblock Homotopy {P}erturbation {M}ethod: {A} {N}ew {N}onlinear {A}nanalytic
  {T}echnique.
\newblock {\em Journal of Applied Mathematics and Computation} 135:73--79.

\bibitem[\protect\astroncite{Hemker}{1990}]{Hemker1990}
{\sc Hemker, P.~W.} 1990.
\newblock A {N}onlinear {M}ultigrid {M}ethod for {O}ne-{D}imensional
  {S}emiconductor {D}evice {S}imulation.
\newblock {\em Journal of Computational and Applied Mathematics} 30:117--126.

\bibitem[\protect\astroncite{Horn and Johnson}{1985}]{Horn1985}
{\sc Horn, R.~A. and Johnson, C.~R.} 1985.
\newblock Matrix {A}nalysis.
\newblock Cambridge University Press.

\bibitem[\protect\astroncite{Le-Gall}{2014}]{LeGall2014}
{\sc Le-Gall, F.} 2014.
\newblock {P}owers of {T}ensors and {F}ast {M}atrix {M}ultiplication.
\newblock {\em {P}roceedings of the 39th {I}nternational {S}ymposium on
  {S}ymbolic and {A}lgebraic {C}omputation ({ISSAC} 2014), ar{X}iv: 1401.7714}
  .

\bibitem[\protect\astroncite{Liao}{1998}]{Liao1998}
{\sc Liao, S.} 1998.
\newblock On the 8th-order {D}rag {C}oefficient of a {S}phere in a {U}niform
  {S}tream: {A} {S}implified {D}escription.
\newblock {\em Communications in Nonlinear Science and Numerical Simulation}
  3(4):256--260.

\bibitem[\protect\astroncite{Liao}{2002}]{Liao2002b}
{\sc Liao, S.} 2002.
\newblock An {A}nalytic {A}pproximation of the {D}rag {C}oefficient for the
  {V}iscous {F}low past a {S}phere.
\newblock {\em International Journal of Non-Linear Mechanics} 37(1):1--18.

\bibitem[\protect\astroncite{Liao}{1992}]{Liao1992}
{\sc Liao, S.~J.} 1992.
\newblock {T}he proposed {H}omotopy {A}nalysis {T}echnique for the solution of
  {N}onlinear problems.
\newblock PhD thesis, Shanghai Jiao Tong University.

\bibitem[\protect\astroncite{Liao}{1997}]{Liao1997}
{\sc Liao, S.~J.} 1997.
\newblock {A} {K}ind of {A}pproximate {S}olution {T}echnique {W}hich {D}oes
  {N}ot {D}epend {U}pon {S}mall {P}arameters ({II}) - {A}n application in
  {F}luid {M}echanics.
\newblock {\em {I}nternational {J}ournal of {N}onlinear {M}echanics}
  32:815--822.

\bibitem[\protect\astroncite{Liao and Campo}{2002}]{Liao2002}
{\sc Liao, S.~J. and Campo, A.} 2002.
\newblock {A}nalytic solutions of the {T}emperature {D}istribution in {B}lasius
  {V}iscous {F}low {P}roblems.
\newblock {\em {J}ournal of {F}luid {M}echanics} 453:411--425.

\bibitem[\protect\astroncite{Lyapunov}{1892}]{Lyapunov1892}
{\sc Lyapunov, A.~M.} 1892.
\newblock The {G}eneral {P}roblem of the {S}tability of {M}otion.
\newblock Taylor \& {F}rancis, {L}ondon.

\bibitem[\protect\astroncite{Motsa et~al.}{2012}]{Motsa2012}
{\sc Motsa, S., Shateyi, S., Marewo, G., and Sibanda, P.} 2012.
\newblock {A}n improved {S}pectral {H}omotopy {A}nalysis {M}ethod for {MHD}
  {F}low in a {S}emi-{P}orous {C}hannel.
\newblock {\em {N}umerical {A}lgorithms} 60:463--481.

\bibitem[\protect\astroncite{Motsa et~al.}{2010}]{Motsa2010}
{\sc Motsa, S., Sibanda, P., and Shateyi, S.} 2010.
\newblock {A} new {S}pectral-{H}omotopy {A}nalysis {M}ethod for {S}olving a
  {N}onlinear {S}econd {O}rder {BVP}.
\newblock {\em {C}ommunications in {N}onlinear {S}cience and {N}umerical
  {S}imulation} 15:2293--2302.

\bibitem[\protect\astroncite{Motsa}{2014}]{Motsa2014}
{\sc Motsa, S.~S.} 2014.
\newblock On the {O}ptimal {A}uxiliary {L}inear {O}perator for the {S}pectral
  {H}omotopy {A}nalysis {M}ethod {S}olution of {N}onlinear {O}rdinary
  {D}ifferential {E}quations.
\newblock {\em Mathematical Problems in Engineering} 2014.

\bibitem[\protect\astroncite{Nik et~al.}{2013}]{Nik2013}
{\sc Nik, H., Effati, S., Motsa, S.~S., and Shirazian, M.} 2013.
\newblock {S}pectral {H}omotopy {A}nalysis {M}ethod and its {C}onvergence for
  solving a class of {N}onlinear {O}ptimal {C}ontrol {P}roblems.
\newblock {\em {N}umerical {A}lgorithms} .

\bibitem[\protect\astroncite{Olver and Townsend}{2013}]{Olver2013}
{\sc Olver, S. and Townsend, A.} 2013.
\newblock {A} {F}ast and {W}ell-{C}onditioned {S}pectral {M}ethod.
\newblock {\em SIAM Review} 55:462--489.

\bibitem[\protect\astroncite{Peyret}{2002}]{Peyret2002}
{\sc Peyret, R.} 2002.
\newblock Spectral {M}ethods for {I}ncompressible {V}iscous {F}low.
\newblock Springer.

\bibitem[\protect\astroncite{R{\o}nquist and Patera}{1987}]{Ronquist1987}
{\sc R{\o}nquist, E. and Patera, A.} 1987.
\newblock Spectral {E}lement {M}ultigrid. {I}: {F}ormulation and {N}umerical
  {R}esults.
\newblock {\em Journal of Scientific Computing} 2(4):389--406.

\bibitem[\protect\astroncite{Rothman}{1991}]{Rothman1991}
{\sc Rothman, E.~E.} 1991.
\newblock Reducing {R}ound-{O}ff {E}rror in {C}hebyshev {P}suedospectral
  {C}omputations.
\newblock Technical report.

\bibitem[\protect\astroncite{Sloane}{2018}]{OEIS}
{\sc Sloane, N. J.~A.} 2018.
\newblock The {O}n-{L}ine {E}ncyclopedia of {I}nteger {S}equences.
\newblock {\em published at http://oeis.org} .

\bibitem[\protect\astroncite{Strikwerda}{2004}]{Strikwerda2004}
{\sc Strikwerda, J.} 2004.
\newblock Finite Difference Schemes and Partial Differential Equations.
\newblock SIAM, 2 edition.

\bibitem[\protect\astroncite{Trefethen and Trummer}{1989}]{Trefethen1989}
{\sc Trefethen, L.~N. and Trummer, M.~R.} 1989.
\newblock An {I}nstability {P}henomenon in {S}pectral {M}ethods.
\newblock {\em SIAM Journal of Numerical Analysis} .

\bibitem[\protect\astroncite{Wang}{2002}]{Wang2002}
{\sc Wang, Z.} 2002.
\newblock Spectral ({F}inite) {V}olume {M}ethod for {C}onservation {L}aws on
  {U}nstructured {G}rids: {B}asic {F}ormulation.
\newblock {\em Journal of Computational Physics} 178:210--251.

\bibitem[\protect\astroncite{Xu et~al.}{2012}]{Xu2012}
{\sc Xu, D.~L., Lin, Z.~L., Liao, S.~J., and Stiassnie, M.} 2012.
\newblock {O}n the {S}teady-{S}tate {F}ully {R}esonant {P}rogressive {W}aves in
  {W}ater of {F}inite {D}epth.
\newblock {\em {J}ournal of {F}luid {M}echanics} 207:1--40.

\bibitem[\protect\astroncite{Yang and Liao}{2006}]{Yang2006}
{\sc Yang, C. and Liao, S.~J.} 2006.
\newblock {O}n the {E}xplicit {P}urely {A}nalytic {S}olution of {V}on {K}arman
  {S}wirling {V}iscous {F}low.
\newblock {\em {C}ommunication in {N}onlinear {S}cience: {N}umerical
  {S}imulations} 11:83--93.

\bibitem[\protect\astroncite{Yang and Lakoba}{2007}]{Yang2007}
{\sc Yang, J. and Lakoba, T.~I.} 2007.
\newblock {A}ccelerated {I}maginary-time {E}volution {M}ethods for the
  {C}omputation of {S}olitary {W}aves.
\newblock {\em ar{X}iv ref:0711.3434v1} .

\bibitem[\protect\astroncite{Yibiao and Xin-Kai}{2016}]{Yibiao2016}
{\sc Yibiao, L. and Xin-Kai, L.} 2016.
\newblock The {C}hebyshev {S}pectral {E}lement {A}pproximation with {E}xact
  {Q}uadratures.
\newblock {\em Journal of Computational and Applied Mathematics} .

\end{thebibliography}







\end{document}